\documentclass[12pt]{article}
\usepackage{amsfonts, amsmath, amsthm, graphics, bbold}

\newcommand{\be}{\begin{equation}}
\newcommand{\eneq}{\end{equation}}

\def\ddx{{{\rm d}\over{\rm d} x}}

\def\thi{{\vartheta_\infty}}
\def\eps{\varepsilon}

\newcommand{\ID}{\mathbb{1}}
\newcommand{\complessi}{\mathbb C}
\newcommand {\interi}{\mathbb Z}

\theoremstyle{definition}
\newtheorem{df}{Definition}

\newtheorem{rmk}[df]{Remark}

\theoremstyle{theorem}
\newtheorem{lm}[df]{Lemma}
\newtheorem{thm}[df]{Theorem}

\newtheorem{prop}[df]{Proposition}
\begin{document}
\title{{\bf The geometry of the classical solutions of the Garnier systems.}}
\author{Marta Mazzocco\thanks{Mathematical
Institute, 24-29 St Giles, Oxford OX1 3LB, UK.}} 
\date{}
\maketitle

\begin{abstract}
Our aim is to find a general approach to the theory of classical 
solutions of the Garnier system in $n$-variables, ${\cal G}_n$, 
based on the Riemann-Hilbert problem and on the geometry of the 
space of isomonodromy deformations. Our approach consists in 
determining the monodromy data of the corresponding Fuchsian system 
that guarantee to have a classical solution of the Garnier system 
${\cal G}_n$. 
This leads to the idea of the {\it reductions of the Garnier systems.}\/ 
We prove that if a solution of the Garnier system ${\cal G}_{n}$ 
is such that the associated Fuchsian system has $l$ monodromy matrices
equal to $\pm\ID$, then it can be reduced classically to a solution
of a the Garnier system with $n-l$ variables ${\cal G}_{n-l}$. When $n$ 
monodromy matrices are equal to $\pm\ID$, we have classical solutions 
of ${\cal G}_n$.
We give also another mechanism to produce classical solutions:
we show that the solutions of the Garnier systems having reducible
monodromy groups can be reduced to the classical solutions found by 
Okamoto and Kimura in terms of Lauricella hypergeometric functions.
In the case of the Garnier system in $1$-variables, i.e. for the  
Painlev\'e VI equation, we prove that all classical non-algebraic 
solutions have either reducible monodromy groups or at least one 
monodromy matrix equal to $\pm\ID$.
\end{abstract}

\section{Introduction.}

%%%%
The $n$--variables {\it Garnier system}\/ ${\cal G}_n$ is the
completely integrable Hamiltonian system \cite{Gar1,Gar2,Ok1}
\begin{equation}
\left\{\begin{array}{cc}
{\partial\nu_j\over\partial u_i}&={\partial\
K_i\over\partial\rho_j}
\qquad i,j=1,\dots,n,\\
{\partial\rho_j\over\partial u_i}&=-{\partial\ K_i\over\partial\nu_j}
\qquad i,j=1,\dots,n,\\
\end{array}
\right.
\label{ga1}\end{equation}
where $u_1,\dots,u_n$ are pairwise distinct
complex variables, and 
\begin{equation}
K_i=-{\Lambda(u_i)\over T'(u_i)}\left[\sum_{k=1}^{n}
{T(\nu_k)\over(\nu_k-u_i)\Lambda'(\nu_k)}
\left\{\rho_k^2-
\sum_{m=1}^{n+2}{\theta_m-\delta_{im}\over\nu_k-u_m}\rho_k
+{\kappa\over\nu_k(\nu_k-1)}\right\}\right]
\label{ham}\end{equation}
with $u_{n+1}=0, u_{n+2}=1$, $\kappa={1\over4}\left\{ \left(
\sum_{m=1}^{n+2}\theta_m-1\right)^2-(\thi+1)^2\right\}$,
$\theta_1,\dots,\theta_{n+2}$, $\thi$ being some constant parameters, 
$\Lambda(u):=\Pi_{k=1}^{n}(u-\nu_k)$ and
$T(u):=\Pi_{k=1}^{n+2}(u-u_k)$.

The Garnier system ${\cal G}_n$ has movable algebraic points but it
can be transformed by a suitable change of variables 
$(\nu_i,\rho_i,u_i)\to(q_i,p_i,s_i)$, $i=1,\dots,n$,
into a Hamiltonian system ${\cal H}_n$ enjoying the 
Painlev\'e property (see \cite{KO}). This means that the solutions 
$q_i(s_1,\dots,s_n),\,p_i(s_1,\dots,s_n)$  may have 
complicated singularities (i.e. branch points, essential singularities etc.) 
only at the {\it critical points}\/ $u_i=u_j$ for $i\neq j$ of 
the equation, the position of which does not depend on the choice of
the particular solution 
(the so-called {\it fixed singularities}\/). All the other singularities, the
position of which depend on the integration constants (the so-called 
{\it movable singularities}\/), are poles.

Observe that in the case of only one variable $n=1$ (and four
parameters $\theta_1$, $\theta_2$, $\theta_3$ and $\thi$) the Garnier
system satisfies the Painlev\'e property also in the variables 
$(\nu,\rho)$ and it coincides with the Painlev\'e sixth equation with 
parameters $\alpha={(\thi-1)^2\over2}$, $\beta=-{\theta_2^2\over 2}$, 
$\gamma={\theta_3^2\over 2}$ and $\delta={1-\theta_1^2\over 2}$. 
The solutions of the Painlev\'e sixth equation are new transcendental 
functions, i.e for generic values of the parameters $\theta_1$, $\theta_2$,
$\theta_3$ and $\thi$, the generic solutions can not be expressed via 
elementary or classical transcendental functions 
(see \cite{Pain, Um1, Wat}). As a consequence one expects that also
the generic solutions of the Garnier system ${\cal G}_n$ with $n>1$
will be some new transcendental functions of several variables.

More precisely, following \cite{Um,Um1,Um2}, a {\it classical 
function}\/ in one variable, $x$ is a function that can be obtained 
from the field of rational functions $\complessi(x)$, by a finite 
iteration of the following {\it admissible operations:}\/
\begin{description}
\item[i)] arithmetic operations $+,-,\times,\div$, 
\item[ii)] quadrature,
\item[iii)] solutions of algebraic equations with classical 
coefficients,
\item[iv)] solution of a linear ordinary differential equations 
with classical coefficients,
\item[v)] substitution into an Abelian functions,
\item[vi)] solution of algebraic differential equations of first 
order with classical coefficients.
\end{description}

For functions of several variables, we give the following

\begin{df}
We say that a function of several variables $f(u_1,\dots,u_n)$
is a {\it classical functions}\/ if for every 
$i=1,\dots,n$, and for every fixed $u_j$, $j\neq i$, the function
$f_i:u_i\to f(u_1,\dots,u_n)$ is classical in $u_i$ in the above sense. 
\label{classic}\end{df}

For some special values of the parameters 
$\theta_1,\dots,\theta_{n+2}$, $\thi$ there are particular solutions 
that can be expressed via classical functions. For example, if 
$\sum_{k=1}^{n+2}\theta_k+1=\pm\thi$, then there are some special
solutions of ${\cal G}_n$ that can be expressed via the Lauricella
hypergeometric functions (see \cite{OK}). In the case of Painlev\'e
sixth equation there are several examples of classical solutions, that
we do not discuss here. We only mention that, in certain cases, the
classical solutions of the Painlev\'e sixth equation have been related 
to the Dynkin diagrams (see \cite{Ok2}) and to the regular polyhedra 
in the Euclidean space (see \cite{DM,M,Hit}). 
Observe that there are some explicit solutions of ${\cal G}_n$ that
are not classical according to Umemura's definition (see 
\cite{Pic,Hit1,KK,DIKZ}). 

The classical solutions of ${\cal G}_n$ have several applications
in the theory of integrable systems and Frobenius manifolds.
Indeed Garnier systems and their limits appear in the context of 
twistor theory as symmetry reductions of several integrable systems.
For example the self--dual Yang--Mills equations with $3$-dimensional 
symmetry groups, the self--dual Einstein equations and various 
generalisations of them have  symmetry reductions to ${\cal G}_n$
(see \cite{MW}). The basic idea of the Painlev\'e test of
integrability of PDEs consists in reducing a given system to a 
differential equation possessing the Painlev\'e property (see for 
example \cite{Con}). In the context of the theory of 
Frobenius manifolds, some Garnier systems with are related to 
the problem of the normal forms of dispersion-less 
integrable systems. In \cite{DZ2} it is shown that this relation remains 
valid for the dispersive corrections. In particular the first dispersive 
correction is expressed by the isomonodromic $\tau$-function. Such an
isomonodromic perspective for the higher order corrections is still
missing. 

Our aim is to find a general approach to the theory of classical 
solutions of the Garnier system based on the Riemann-Hilbert problem 
and on the geometry of the space of isomonodromy deformations (see
\cite{FN, JMU}). 
In fact systems where first introduced by Garnier as isomonodromic
deformations equations of a certain second order Fuchsian differential
equation, or, equivalently of an auxiliary two-dimensional Fuchsian 
system with $n+3$ poles $u_1,\dots,u_{n+2},\infty$ 
(see \cite{Gar1,Gar2,Ok1,IKSY}).

Our approach consists in determining the monodromy data of the
Fuchsian equation that guarantee to have a classical solution. This
methods enables us to fully reproduce the results obtained in \cite{OK} 
in the framework of the symmetries of Hamiltonian systems. Moreover
our approach suggests the {\it reductions of the Garnier systems.}\/ 
In fact we show that all solutions of the Garnier system 
${\cal G}_n$ can be seen as solutions of the Garnier system 
${\cal G}_{n+1}$ but there are solutions of the Garnier system 
${\cal G}_{n+1}$ that can not be reduced to solutions of the Garnier 
system ${\cal G}_{n}$ via classical operations. In this sense we expect 
that the generic solutions of the Garnier system ${\cal G}_{n+1}$ 
should be {\it more transcendental}\/ than the ones of ${\cal G}_n$.

More precisely, the Garnier system ${\cal G}_n$ is represented as 
the equation of isomonodromy deformation of the two-dimensional
auxiliary Fuchsian system with $2\times 2$ residue matrices 
${\cal A}_j$ independent on $\lambda$:
\begin{equation}
{{\rm d}\over{\rm d}\lambda} \Phi=\left(
\sum_{k=1}^{n+2}{{\cal A}_k\over \lambda-u_k} \right)\Phi,
\label{N1in}
\end{equation}
$u_1,\dots,u_{n+2}$ being pairwise distinct complex numbers. 
The residue matrices ${\cal A}_j$ satisfy the following conditions:
$$
{\rm eigen}\left({\cal A}_j\right)=\pm{\theta_j\over2}
\quad\hbox{and}\quad
-\sum_{k=1}^{n+2}{\cal A}_k={\cal A}_\infty,
$$
where
$$
\hbox{for}\,\theta_\infty\neq 0,\quad
{\cal A}_\infty:={1\over2}\left(
\begin{array}{cc}
\theta_\infty & 0\\ 0 & -\theta_\infty\\
\end{array}\right),
\quad\hbox{and for}\,\thi=0,
\quad  
{\cal A}_\infty:=\left(
\begin{array}{cc}0& 1\\ 0 & 0\\\end{array}\right),
$$
$\theta_j$, $j=1,\dots,n+2,\infty$ are the ones appearing in (\ref{ham}).

Generically (see \cite{Bol}) the monodromy matrices 
${\cal M}_1,\dots,{\cal M}_{n+2}$ 
of (\ref{N1in}) remain constant if and only if the residue matrices 
${\cal A}_1,\dots,{\cal A}_{n+2}$ are solutions of the Schlesinger
equations (see \cite{Sch}):
\begin{equation}
{\partial\over\partial u_j} {\cal A}_i= 
{[ {\cal A}_i, {\cal A}_j]\over u_i-u_j},\quad
{\partial\over\partial u_i} {\cal A}_i= 
-\sum_{j\neq i}{[ {\cal A}_i, {\cal A}_j]\over u_i-u_j}. 
\label{schleq}\end{equation}
When we deal with $2\times 2$ matrices, the above equations reduce 
to the Garnier system ${\cal G}_n$ (see \cite{Gar1,Gar2,IKSY}).

We show that when the  monodromy group 
$\langle{\cal M}_1,\dots,{\cal M}_{n+2}\rangle$ is 
{\it $l$-smaller,}\/ i.e. it is such that $l$ generating monodromy 
matrices are equal to $\pm\ID$, then 
the corresponding solutions of the Garnier system ${\cal G}_n$
{\it reduce}\/ to solutions of the Garnier system ${\cal G}_{n-l}$.
If the monodromy group is isomorphic to the monodromy group of the 
hypergeometric equation (i.e. is $n$-smaller), we have classical 
solutions. More precisely we prove the following

\begin{thm}
If there exists a solution $(\nu_1,\dots,\nu_n,\rho_1,\dots,\rho_n)$
of the Garnier system ${\cal G}_{n}$ such that the corresponding 
Fuchsian system of the form (\ref{N1in}) with matrices 
$({\cal A}_1,\dots,{\cal A}_{n+2})$ is $l$-smaller, then there is 
a $l$-parameters 
family of such solutions and there exists a solution $(\hat\nu_1,
\dots,\hat\nu_{n-l},\hat\rho_1,\dots,\hat\rho_{n-l})$
of the Garnier system ${\cal G}_{n-l}$ such that the corresponding
Fuchsian system has monodromy group generated by those matrices 
${\cal M}_1,\dots,{\cal M}_{n+2}$ that are not equal to $\pm\ID$. 
Moreover the obtained $l$-parameter family of solutions of the 
Garnier system ${\cal G}_{n}$ depends classically on 
$(\hat\nu_1,\dots,\hat\nu_{n-l},\hat\rho_1,
\dots,\hat\rho_{n-l},u_1,\dots,u_{n+2})$. If $l=n$ then we obtain a
$n$-parameter family of classical functions of $(u_1,\dots,u_{n+2})$.
\label{mainthmsm}\end{thm}

\begin{rmk}
Observe that the existence on an $l$-parameters family of solutions
with the same monodromy matrices is related to the fact that 
the theorem of uniqueness of a Fuchsian system with prescribed poles and 
monodromy fails when one or more of the monodromy matrices is a
multiple of the identity (see section 2.1). The dependence of the
Fuchsian system on the $l$-parameters is an example of non-Schlesinger
isomonodromic deformation (see \cite{Bol}).
\label{rm0.1.0}\end{rmk}

\begin{rmk}
In the case when the monodromy at infinity ${\cal M}_\infty$ is 
equal to $\pm\ID$, it is possible to prove an analogous result to 
Theorem \ref{mainthmsm}. The statement is a little more delicate 
and can be found in section 3.2.
\label{rm0.0.0}\end{rmk}

\begin{rmk}
The fact that $(\nu_1,\dots,\nu_n,\rho_1,\dots,\rho_n)$ are classical 
functions of the arguments 
$(\hat\nu_1,\dots,\hat\nu_{n-l},\hat\rho_1,
\dots,\hat\rho_{n-l}, u_1,\dots,u_{n+2})$ does not necessarily imply 
that they are classical functions of the variables $(u_1,\dots,u_{n+2})$.  
\label{rm0.0}\end{rmk}

We give also another mechanism to produce classical solutions:
we show that the solutions of the Garnier systems having reducible
monodromy groups can be reduced to the classical solutions found by 
Okamoto and Kimura in terms of Lauricella hypergeometric functions.

\begin{thm}
If there exists a solution $(\nu_1,\dots,\nu_n,\rho_1,\dots,\rho_n)$
of the Garnier system ${\cal G}_{n}$ with parameters 
$\theta_1,\dots,\theta_{n+2},\thi$ such that the corresponding 
Fuchsian system has a reducible monodromy group then it is a 
classical solution. For such parameters 
$\theta_1,\dots,\theta_{n+2},\thi$ there exists a $n$-parameter family
of classical solutions of the Garnier system ${\cal G}_{n}$ all having
reducible monodromy groups.
\label{mainthmred}\end{thm}

In the case of the Painlev\'e VI equation, we can prove that

\begin{thm}
All classical non-algebraic solutions of the Painlev\'e VI equation  
have either reducible or smaller monodromy groups.
\label{thmio}\end{thm}

It is then natural to believe that all solutions of the Garnier system 
with {\it generic monodromy groups,}\/ i.e. monodromy groups that 
are non-reducible and non-smaller, are either algebraic or non
classical. We cannot prove such a result yet, but this paper gives 
arguments that make such a belief stronger.

Observe that when the matrices ${\cal A}_1,\dots,{\cal A}_{n+2}$
are $m\times m$, the generic isomonodromic deformations equations are
still given by the Schlesinger equations (\ref{schleq}), i.e. by 
the so-called {\it Schelsinger systems.}\/ These systems have very 
rich analytical and geometric structures. For example they admit a 
Hamiltonian formulation whose quantisation can be regarded as 
Knizhnik--Zamolodchikov system and its generalisations (see \cite{Res}).
We postpone to another paper the study of the reductions and
classical solutions of the Schlesinger systems (see \cite{DM1}).

\vskip 0.5 cm
\noindent{\bf Acknowledgements.} The author is grateful to
B. Dubrovin and M. Jimbo who inspired this paper and is indebted to
N. Hitchin, who constantly addressed this work and gave lots of 
suggestions. The author was supported by an EPSRC research 
assistantship and partially by SISSA, International School of Advanced
Studies, Trieste and by Worcester College, Oxford.

%%%%%%%%%%%%%%%%%%%%%%%%%%%%

%%%%%%%% SEZIONE 2

%%%%%%%%%%%%%%%%%%%%%%%%%

\section{Isomonodromic deformations equations and Garnier systems}

In this section we recall how to represent the Garnier systems as 
isomonodromic deformation equations of a $2\times2$ Fuchsian system. 
The results of this section are standard 
(see \cite{Gar1,Gar2,Ok1,KO,IKSY}), we just recall them here in order
to fix the notations and for self-consistency.

\subsection{Monodromy data of Fuchsian systems}

Consider the Fuchsian system with $n+3$ regular singularities at 
$u_1,\dots$, $u_{n+2}$, $u_{n+3}=\infty$:
\begin{equation}
{{\rm d}\over{\rm d}\lambda} \Phi=
\sum_{k=1}^{n+2}{{\cal A}_k\over \lambda-u_k}\Phi,
\qquad\qquad \lambda\in\overline\complessi\backslash\{u_1,\dots,u_{n+3}\}
\label{N1}\end{equation}
${\cal A}_j$ being $2\times 2$ matrices independent on $\lambda$, and 
$u_i\neq u_j$ for $i\neq j$, $i,j=1,\dots,n+2$. Take some parameters
$\theta_1,\dots,\theta_{n+2},\thi$ and assume that 
the matrices ${\cal A}_j$ satisfy the following conditions:
\begin{equation}
{\rm eigen}\left({\cal A}_j\right)=\pm{\theta_j\over2},
\quad\hbox{and}\quad
-\sum_{k=1}^{n+2}{\cal A}_k={\cal A}_\infty,
\label{N1.3}
\end{equation}
where 
$$
{\cal A}_\infty:={1\over2}\left(
\begin{array}{cc}
\theta_\infty & 0\\ 0 & -\theta_\infty\\
\end{array}\right),
$$
for some $\theta_\infty\neq 0$.

The solution $\Phi(\lambda)$ of the system (\ref{N1}) is a multi-valued 
analytic function in the punctured Riemann sphere 
$\complessi\backslash\{u_1,\dots u_{n+2}\}$, and its multivaluedness 
is described by the so-called {\it monodromy matrices.}\/  Let us briefly 
recall the definition of the monodromy matrices of the Fuchsian system 
(\ref{N1}). First, fix a basis $\gamma_1,\dots,\gamma_{n+2}$ of loops
in $\pi_1\left(\overline\complessi\backslash\{u_1,\dots u_{n+3}\},
\infty\right)$, and a fundamental matrix for the system (\ref{N1}).
To fix the basis of the loops, one has to perform some cuts between the 
singularities, i.e. $n+1$ segments $\pi_1\dots\pi_{n+2}$ between 
$u_{n+3}=\infty$ and each $u_j$, $j=1,\dots,n+2$. The segments
$\pi_{1,\dots,n+2}$ are taken along the same direction $\eta$ and ordered 
according to the order of the points $u_1,\dots,u_{n+2}$. 
% I fix the order of
% the points $u_1,u_2,u_3$ to be {\it lexicographical,}\/ i.e. 
% $\RE u_k\leq \RE u_{k+1}$, and if $\RE u_k=\RE u_{k+1}$ then 
% $\IM u_k<\IM u_{k+1}$.
Take $\gamma_j$ to be a simple closed curve starting and finishing at 
infinity, going around $u_j$ in positive direction ($\gamma_j$ is oriented 
counter-clockwise, $u_j$ lies inside, while the other singular points lie 
outside) and not crossing the cuts $\pi_i$. Near $\infty$, every 
loop $\gamma_j$ is close to the cut $\pi_j$ as in Figure 1.

%%%%

%%%%%% figura F1. %%%%%%%%%%%%

\begin{figure}
\includegraphics{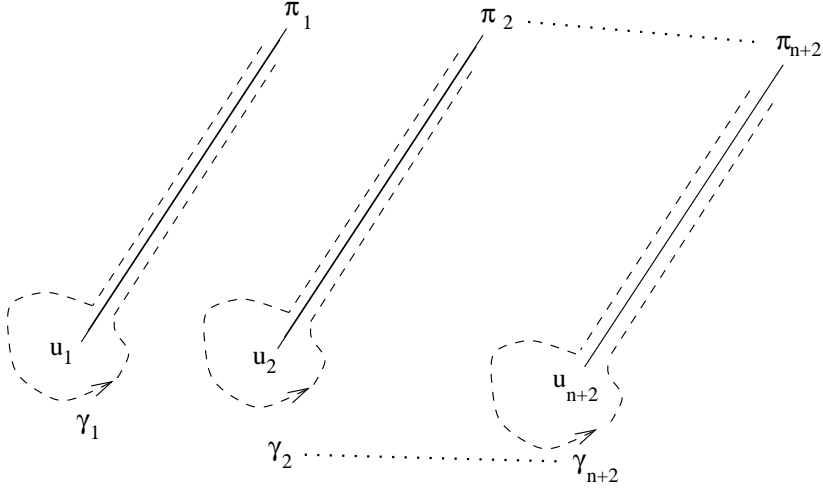}
\caption{The basis of loops $\gamma_1,\dots,\gamma_{n+2}$.}
\end{figure}

The fundamental matrix of the system is given by
the following

%%%%%%%%%%% proposition

\begin{prop} There exists a fundamental matrix of the system 
(\ref{N1}) of the form
\begin{equation}
\Phi_\infty(\lambda)=\left(\ID +{\cal O}({1\over \lambda})\right) 
\lambda^{-{\cal A}_\infty}\lambda^{-{\cal R}_\infty},\quad\hbox{as}\quad 
\lambda\rightarrow\infty,
\label{M3}
\end{equation}
where the matrix ${\cal R}_\infty$ has zero diagonal entries and
\begin{equation}\begin{array}{l}
{\cal R}_{\infty_{ij}} = 0 \qquad \hbox{if}\quad
\thi\not\in\interi,\\
{\cal R}_{\infty_{12}} =
-\sum_{k=1}^{n+2}\left({\cal A}_k\right)_{12} u_k^p-
\sum_{l=1}^{p-1}\left(G^{(p-l)}
\sum_{k=1}^{n+2} {\cal A}_k u_k^l\right)_{12},
\quad \hbox{if}\,
\thi=p\in\interi_+,\\ 
{\cal R}_{\infty_{21}} =
-\sum_{k=1}^{n+2}\left({\cal A}_k\right)_{21} u_k^p-
\sum_{l=1}^{p-1}\left(G^{(p-l)}
\sum_{k=1}^{n+2} {\cal A}_k u_k^l\right)_{21},
\quad \hbox{if}\,
\thi=-p\in\interi_-,\\
\end{array}
\end{equation}
where $G^{(0)}=\ID$ and for $l=1,2,\cdots p-1$, $G^{(l)}$ are uniquely 
determined by
\begin{equation}
[{\cal A}_\infty,G^{(l)}]=l G^{(l)} -
\sum_{k=1}^{n+2}{\cal A}_k u_k^l-
\sum_{s=1}^{l-1} G^{(l-s)}
\sum_{k=1}^{n+2} {\cal A}_k u_k^s
\end{equation}
and $\lambda^{-{\cal R}_\infty}:=e^{-{\cal R}_\infty\log\lambda}$, with 
the choice of the principal branch of the logarithm with the
branch-cut along the common direction $\eta$ of the cuts 
$\pi_1,\dots,\pi_{n+2}$. Such a fundamental matrix 
$\Phi_\infty(\lambda)$ is uniquely determined up to 
\begin{equation}
\Phi_\infty(\lambda)\to \Phi_\infty(\lambda)L_\infty,
\label{amb}
\end{equation}
where $L_\infty$ is any constant invertible matrix such that 
\begin{equation}
\lambda^{-{\cal R}_\infty}\lambda^{-{\cal A}_\infty}L_\infty
\lambda^{{\cal A}_\infty}\lambda^{{\cal R}_\infty}=
\ID + \sum_{k=1}^N {L_\infty^{(k)}\over\lambda^k},
\label{amb1}
\end{equation}
for some $L_\infty^{(1)},\cdots,L_\infty^{(N)}$ constant matrices.
\label{prop1}\end{prop}

\noindent The proof can be found in \cite{Dub1}. 
\vskip 0.2 cm 

%%%%%%%% monodromy

The fundamental matrix $\Phi_\infty$ can be analytically 
continued to an analytic function on the universal covering of 
$\overline\complessi\backslash\{u_1,\dots,u_{n+3}\}$. For any element 
$\gamma\in\pi_1
\left(\overline\complessi\backslash\{u_1,\dots,u_{n+3}\},\infty\right)$
denote by $\gamma[\Phi_\infty(\lambda)]$ the result of the analytic 
continuation of $\Phi_\infty(\lambda)$ along the loop $\gamma$. Since 
$\gamma[\Phi_\infty(\lambda)]$ and $\Phi_\infty(\lambda)$ are two fundamental 
matrices in the neighbourhood of infinity, they are related by the 
following relation:
\begin{equation}
\gamma[\Phi_\infty(\lambda)]= \Phi_\infty(\lambda) {\cal M}_\gamma
\end{equation}
for some constant invertible $2\times 2$ matrix ${\cal M}_\gamma$
depending only 
on the homotopy class of $\gamma$. Particularly, the matrix 
${\cal M}_\infty:={\cal M}_{\gamma_\infty}$, $\gamma_\infty$ being a
simple loop around 
infinity in the clock-wise direction, is given by:
\begin{equation}
{\cal M}_\infty=\exp(2\pi i {\cal A}_\infty)
\exp(2\pi i{{\cal R}_\infty}).
\label{N9}
\end{equation}
The resulting {\it monodromy representation}\/ is an anti-homomorphism:
\begin{equation}
\begin{array}{ccc}
\pi_1\left(\overline\complessi\backslash
\{u_1,,\dots,u_{n+2},\infty\},\infty\right)
&\rightarrow & SL_2(\complessi)\\
\gamma&\mapsto & {\cal M}_\gamma\\
\end{array}
\label{N2.2}
\end{equation}
\begin{equation}
{\cal M}_{\gamma\tilde\gamma}= {\cal M}_{\tilde\gamma}  
{\cal M}_{\gamma}. \label{N2}
\end{equation}
The images ${\cal M}_j:={\cal M}_{\gamma_j}$ of the generators 
$\gamma_j$, $j=1,\dots,n+2$ of 
the fundamental group, are called {\it the monodromy matrices}\/ of 
the Fuchsian system (\ref{N1}). They generate the 
{\it monodromy group of the system,}\/ i.e. the image of the representation 
(\ref{N2.2}). Since the loop $(\gamma_1 \cdots \gamma_{n+2})^{-1}$ is
homotopic to $\gamma_\infty$, the following relation between the
generators holds:
\begin{equation}
{\cal M}_\infty {\cal M}_{n+2} \cdots {\cal M}_1=\ID.
\label{N6}
\end{equation}
Observe that if we fix another fundamental matrix 
$\Phi_\infty'=\Phi_\infty L_\infty$ in the equivalence class defined 
by (\ref{amb}), the monodromy matrices ${\cal M}_\gamma'$ with
respect to the new fundamental matrix $\Phi_\infty'$ are related to
the old ones by
\begin{equation}
{\cal M}_j'= L_\infty^{-1}{\cal M}_jL_\infty,\quad
j=1,\dots,n+2.
\end{equation}
Thus, given the Fuchsian system (\ref{N1}), with the constraints
(\ref{N1.3}), and the basis of loops as in Figure 1, the monodromy 
matrices ${\cal M}_j$, $j=1,\dots,n+2$ are uniquely defined up to 
the ambiguity
\begin{equation}
({\cal M}_1, \dots, {\cal M}_{n+2})\sim
(L_\infty^{-1}{\cal M}_1L_\infty, \dots,
L_\infty^{-1}{\cal M}_{n+2}L_\infty),
\label{amb2}
\end{equation}
where $L_\infty$ is given by (\ref{amb1}). Observe that ${\cal M}_\infty$ 
is invariant w.r.t. (\ref{amb2}).

We recall the definition of the {\it connection matrices.}\/ Near the poles 
$u_k$, the fundamental matrices $\Phi_k(\lambda)$ of the 
system (\ref{N1}), are given by the following

%%%%%%%%%% proposition

\begin{prop}
There exists a fundamental matrix of the system 
(\ref{N1}) of the form
\begin{equation}
\Phi_k(\lambda)={\cal G}_k\left(\ID +{\cal O}(\lambda-u_k)\right) 
(\lambda-u_k)^{J_k}(\lambda-u_k)^{{\cal R}_k},\quad\hbox{as}\quad 
\lambda\rightarrow u_k,
\label{N6.1}
\end{equation}
where $J_k$ is the Jordan normal form of ${\cal A}_k$, with
eigenvalues $\pm{\theta_k\over2}$, ${\cal G}_k$ is 
defined by ${\cal A}_k={\cal G}_k J_k {\cal G}_k^{-1}$, and the matrix
${\cal R}_k$ has zero diagonal elements and off-diagonal ones given by
\begin{equation}\begin{array}{lc}
{\cal R}_{k_{ij}} = 0 &\quad \hbox{if}\quad
\theta_k\not\in\interi,\\
{\cal R}_{k_{12}} = \left ( \hat{\cal A}_{k,p}+
\sum_{l=1}^{n-1}\hat {\cal A}_{k,n-l} G_k^{(l)}\right)_{12}
&\quad \hbox{if}\quad
\theta_k=p\in\interi_+,\\
{\cal R}_{k_{21}} = \left ( \hat{\cal A}_{k,p}+
\sum_{l=1}^{n-1}\hat {\cal A}_{k,n-l} G_k^{(l)}\right)_{21}
&\quad \hbox{if}\quad
\theta_k=-p\in\interi_-,\\ \end{array}
\end{equation}
where 
\begin{equation}
\hat {\cal A}_{k,l}= {{\cal G}}_k \sum_{j\neq k} (-1)^{l-1} 
{{\cal A}_j\over(u_k-u_j)^{l}}{{\cal G}}_k^{-1}
\end{equation}
and  $G_k^{(0)}=\ID$ and $G_k^{(l)}$, for $l=1,2,\cdots n-1$,  are 
uniquely determined by
\begin{equation}
0= [J_k, G_k^{(l)}] -l G_k^{(l)}+
\sum_{m=1}^{l-1}\hat {\cal A}_{k,l-m} G_k^{(m)} + \hat {\cal A}_{k,l}.
\end{equation}
The choice of the branch of $\log(z-u_k)$ needed in the definition of
$(\lambda-u_k)^{J_k}$ and $(\lambda-u_k)^{{\cal R}_k}$ is the same as
in Proposition \ref{prop1}. The fundamental matrix $\Phi_k(\lambda)$ 
is uniquely determined up to the ambiguity:
\begin{equation}
\Phi_k(\lambda)\mapsto \Phi_k(\lambda)L_k
\end{equation}
where $L_k$ is any constant invertible matrix such that 
\begin{equation}
(\lambda-u_k)^{J_k}(\lambda-u_k)^{{\cal R}_k}L_k
(\lambda-u_k)^{-J_k}(\lambda-u_k)^{-{\cal R}_k}=
\sum_{j=0}^N {L_k^{(j)}(\lambda-u_k)^j},\label{amb5}
\end{equation}
for $L_k^{(0)}={\cal G}_k$ and for some $L_k^{(1)},\cdots,L_k^{(N)}$ constant
matrices.
\label{prop2}
\end{prop}

\noindent The proof can be found in  \cite{Dub1}.
\vskip 0.3 cm

%%%%%%%%%%%% connection matrices

Continuing the solution $\Phi_\infty(\lambda)$ to a neighbourhood of $u_j$, 
along, say, the right-hand-side of the cut $\pi_j$, one obtains another 
fundamental matrix around $u_j$, that must be related to $\Phi_j(\lambda)$ by:
\begin{equation}
\Phi_\infty(\lambda)=\Phi_j(\lambda){\cal C}_j,\label{N3}
\end{equation}
for some invertible matrix ${\cal C}_j$. The matrices 
${\cal C}_1,\dots,{\cal C}_{n+2}$  are called 
{\it connection matrices,}\/ and they are defined by (\ref{N3}) up to
the ambiguity ${\cal C}_j\to{\cal C}_jL_\infty$ due to (\ref{amb}).
The connection matrices are related to the monodromy matrices as 
follows:
\begin{equation}
{\cal M}_j={\cal C}_j^{-1} \exp\left(2\pi i J_j\right)
\exp\left( {\cal R}_j\right) 
{\cal C}_j,\qquad j=1,\dots,n+2.
\label{N4}
\end{equation}
Thanks to the above relation it follows that
\begin{equation}
{\rm eigen}({\cal M}_j)=\exp(\pm\pi i\theta_j).
\label{N5}
\end{equation}
%%

%%%%%% lemma

\begin{lm}
Given $n+2$ matrices 
${\cal M}_1,\dots, {\cal M}_{n+2}$, none of which is equal to 
$\pm\ID$, satisfying the relations (\ref{N6}) and (\ref{N5}), then
\item{i)} there exist $n+2$ matrices ${\cal C}_1,\dots,{\cal C}_{n+2}$ 
satisfying the (\ref{N4}). Moreover the matrices 
${\cal C}_1,\dots,{\cal C}_{n+2}$ are uniquely determined by 
the matrices ${\cal M}_1,\dots,{\cal M}_{n+2}$, up to the ambiguity 
${\cal C}_j\mapsto L_j^{-1} {\cal C}_j$, where $L_j$ is any invertible
matrix satisfying (\ref{amb5}).
\item{ii)} If the matrices ${\cal M}_1,\dots,{\cal M}_{n+2}$ 
are the monodromy matrices of a Fuchsian system of the form
(\ref{N1}), then any $(n+2)$-ple ${\cal C}_1,\dots,{\cal C}_{n+2}$
satisfying (\ref{N4}) can be realized as the connection matrices of 
the Fuchsian system itself. 
\label{lm2.5}
\end{lm}

\noindent Proof. i) If ${\cal M}_j$ is diagonalisable, then ${\cal C}_j$
is its diagonalising matrix. If ${\cal M}_j$ is not diagonalisable, it
can be reduced to the Jordan normal form. We can always choose the 
Jordan normal form in such a way that the off-diagonal elements are
all equal to $2\pi i$. Then ${\cal C}_j$ is the matrix reducing 
${\cal M}_j$ to the Jordan normal form chosen in this way.
Two matrices ${\cal C}_j$ and ${\cal C}'_j$ give the same matrix 
${\cal M}_j$ if and only if ${\cal C}_j^{-1} {\cal C}'_j$ commutes
with $\exp\left(2\pi i J_j\right)\exp\left({\cal R}_j\right)$, 
namely if and only if they are related by 
${\cal C}_j = L_j^{-1} {\cal C}'_j$.

ii) Now assume that ${\cal C}'_1,\dots, {\cal C}'_{n+2}$ 
are the connection matrices of a Fuchsian system of the form 
(\ref{N1}), with monodromy matrices ${\cal M}_1,\dots,{\cal M}_{n+2}$ 
with respect to a fixed fundamental matrix $\Phi_\infty$. Id est, 
$\Phi_\infty(\lambda)=\Phi_j'(\lambda){\cal C}'_j$, $j=1,\dots,n+2$, for
some choice of the solutions $\Phi_1',\dots,\Phi_{n+2}'$ of the
form (\ref{N6.1}). One has for each $j=1,\dots,n+2$
\begin{equation}
{\cal M}_j=
({\cal C}'_j)^{-1} \exp\left(2\pi i(J_j+{\cal R}_j)\right){\cal C}'_j=
{\cal C}_j^{-1}\exp\left(2\pi i(J_j+{\cal R}_j)\right){\cal C}_j.
\end{equation}
So the matrices $L_j={\cal C}'_j {\cal C}_j^{-1}$ commute with 
$J_j+{\cal R}_j$ and ${\cal C}_1,\dots,{\cal C}_{n+2}$ are the
connection matrices with respect to the new 
local solutions $\Phi_j(\lambda)=\Phi_j'(\lambda)L_j$.
{\hfill $\bigtriangleup$}

%%%%% Monodromy data

\begin{df}
The {\it Monodromy data}\/ of the Fuchsian system (\ref{N1}) are 
\begin{equation}
\left\{({\cal M}_1,\dots,{\cal M}_{n+2})
\slash_\sim,\,
{\cal R}_1,\dots,{\cal R}_{n+2}\right\},
\end{equation}
where $\sim$ is the equivalence relation defined by (\ref{amb2}).
\label{df2.6}
\end{df}

\begin{rmk}
For non-resonant ${\cal A}_j$, i.e. for $\theta_j\not\in\interi$, 
the correspondent ${\cal R}_j$ matrix is zero by definition and does 
not appear in the set of the monodromy data.
\label{rmk2.7}
\end{rmk}

%%%%%%%% uniqueness lemma

\begin{lm}
Two Fuchsian systems (\ref{N1}) with the same poles
$u_1,\dots,u_{n+2}$, the same exponents $\theta_j$,
$j=1,\dots,n+2$ and the same ${\cal A}_\infty$, having all monodromy
matrices ${\cal M}_1,\dots,{\cal M}_{n+2}$ different from $\pm\ID$, 
coincide if and only if they have the same monodromy data with respect 
to the same basis of the loops $\gamma_1,\dots,\gamma_{n+2}$ given 
in Figure 1. 
\label{lm2.8}
\end{lm}

\noindent Proof. 
Let $\Phi_\infty^{(1)}(\lambda)$ and $\Phi_\infty^{(2)}(\lambda)$ be
the fundamental matrices of the form (\ref{M3}) of the two Fuchsian 
systems. Consider the following matrix:
\begin{equation}
Y(\lambda):= \Phi_\infty^{(2)}(\lambda)\Phi_\infty^{(1)}(\lambda)^{-1}.
\end{equation}
$Y(\lambda)$ is an analytic function around infinity: 
\begin{equation}
Y(\lambda)=1+{\cal O}\left({1\over\lambda}\right),\quad\hbox{as}\, 
\lambda\rightarrow\infty.
\end{equation}
Since the monodromy matrices coincide, $Y(\lambda)$ is a single valued 
function on $\overline\complessi\backslash\{u_1,\dots,u_{n+2}\}$. 
We prove that $Y(\lambda)$ is analytic also at the points $u_j$. In fact having
fixed the monodromy data, we can choose the fundamental matrices 
$\Phi_j^{(1)}(\lambda)$ and $\Phi_j^{(2)}(\lambda)$ of the form
(\ref{N6.1}) with the same exponents ${\cal R}_j\neq0$ and, due to 
Lemma \ref{lm2.5}, in such a way that
\begin{equation}
\Phi_\infty^{(1),(2)}(\lambda)= 
\Phi_j^{(1),(2)}(\lambda){\cal C}_j \quad j=1,2,3.
\end{equation}
with the same connection matrices ${\cal C}_j$. Then near the point
$u_j$, $Y(\lambda)$ is analytic:
\begin{equation}
Y(\lambda)=G_j^{(2)}\left(\ID+{\cal O}(\lambda-u_j) \right)
\left[G_j^{(1)}\left(\ID+{\cal O}(\lambda-u_j) \right)\right]^{-1}.
\end{equation}
This proves that $Y(\lambda)$ is an analytic function on all 
$\overline\complessi$ and then, by the Liouville theorem
$Y(\lambda)=\ID$, and the two Fuchsian systems coincide. 
{\hfill $\bigtriangleup$}

%%%%%%%%%% isomonodromic deformations

\subsection{Isomonodromic deformation equations.}

The theory of the deformations the poles of the Fuchsian system keeping the 
monodromy fixed is described by the following two results:

\begin{thm}
Let $\left\{({\cal M}_1,\dots,{\cal M}_{n+2})
\slash_\sim, \, {\cal R}_1,\dots,{\cal R}_{n+2}\right\}$, be
monodromy data of the Fuchsian system:
\begin{equation}
{{\rm d}\over{\rm d}\lambda} \Phi^0=
\sum_{k=1}^{n+2}{{\cal A}^0_k\over \lambda-u^0_k}\Phi^0,
\label{N7}
\end{equation}
of the above form (\ref{N1.3}), with pairwise distinct poles $u_j^0$, 
and with respect to some basis $\gamma_1,\dots,\gamma_{n+2}$ of the loops in 
$\pi_1\left(\overline\complessi\backslash\{u^0_1,\dots,u^0_{n+3}\},
\infty\right)$. If none of the monodromy matrices is equal to
$\pm\ID$, then there exists an open neighbourhood $U\subset\complessi^{n+2}$ 
of the point $u^0=(u^0_1,\dots,u^0_{n+2})$ such that, for any 
$u= (u_1,\dots,u_{n+2})\in U$, there exists a unique $(n+2)$-ple
${\cal A}_1(u),\dots, {\cal A}_{n+2}(u)$ of analytic matrix
valued functions such that:
\begin{equation}
{\cal A}_j(u^0)= {\cal A}_j^0,\quad i=1,\dots,n+2,
\end{equation}
and the monodromy matrices of the Fuchsian system
\begin{equation}
{{\rm d}\over{\rm d}\lambda} \Phi= \sum_{k=1}^{n+2}
{{\cal A}_k(u)\over \lambda-u_k}\Phi, \label{N8}
\end{equation}
with respect to the same basis\footnote{Observe that the 
basis $\gamma_1,\dots,\gamma_{n+2}$ of 
$\pi_1\left(\overline\complessi\backslash\{u_1,\dots,u_{n+3}\},
\infty\right)$ varies continuously with small variations of 
$u_1,\dots,u_{n+2}$. This new basis is homotopic to the initial one, 
so one can identify them.} $\gamma_1,\dots,\gamma_{n+2}$ of the loops, 
coincide with the given ${\cal M}_1,\dots,{\cal M}_{n+2}$.
The matrices ${\cal A}_j(u)$ are the solutions of the Cauchy problem with 
the initial data $ {\cal A}_j^0$ for the following Schlesinger equations:
\begin{equation}
{\partial\over\partial u_j} {\cal A}_i= 
{[ {\cal A}_i, {\cal A}_j]\over u_i-u_j},\quad
{\partial\over\partial u_i} {\cal A}_i= 
-\sum_{j\neq i}{[ {\cal A}_i, {\cal A}_j]\over u_i-u_j}. 
\label{N10}
\end{equation}
The solution $\Phi_\infty^0(\lambda)$ of (\ref{N7}) of the form 
(\ref{M3}) can be uniquely continued, for $\lambda\neq u_i$ $i=1,\dots,n+2$, 
to an analytic function $\Phi_\infty(\lambda,u),\quad u\in U$, such that
\begin{equation}
\Phi_\infty(\lambda,u^0)=\Phi_\infty^0(\lambda).
\end{equation}
This continuation is the local solution of the Cauchy problem with 
the initial data $\Phi_\infty^0(\lambda)$ for the following system 
that is compatible to the system (\ref{N8}):
\begin{equation}
{\partial\over\partial u_i} \Phi = -{ {\cal A}_i(u)\over \lambda-u_i} \Phi.
\end{equation}
Moreover the functions $ {\cal A}_i(u)$ and $\Phi_\infty(\lambda,u)$ can be 
continued analytically to global meromorphic functions on the universal 
coverings of
\begin{equation}
\complessi^3\backslash\{diags\}:=
\left\{(u_1,\dots,u_{n+2})\in\complessi^{n+2}\,|\,u_i\neq u_j\,\hbox{for}\, 
i\neq j\right\},
\end{equation}
and
\begin{equation}
\left\{(\lambda,u_1,\dots,u_{n+2})\in\complessi^4\,
|\,u_i\neq u_j\quad\hbox{for}\quad i\neq j\,\hbox{and}\,
\lambda\neq u_i,\, i=1,\dots,n+2\right\},
\end{equation}
respectively.
\label{thm2.10}
\end{thm}

The proof of this theorem can be found, for example, in \cite{Mal,Mi,Sib}.
As shown in \cite{Bol}, there are some non-generic situations in which 
the Schlesinger equations (\ref{N10}) do not describe all
isomonodromic deformations of the system (\ref{N1}). In this paper we
deal only with Schlesinger isomonodromic deformations.

\begin{df}
We call the system of differential equations (\ref{N10}) in $n$
variables $u_1,\dots,u_n$, for the $m\times m$ matrices 
${\cal A}_1,\dots,{\cal A}_{n+2}$, {\it Schlesinger system} 
${\cal S}_{(n,m)}$. 
\label{df22}
\end{df}

In this paper we deal with ${\cal S}_{(n,2)}$, i.e. $2\times 2$ 
matrices. The case of higher dimensional matrices is postponed to
another paper (see \cite{DM1}).

The solvability of the Schlesinger systems with given monodromy
matrices is still an open problem. Existence can be prove in generic
cases, as in the following

\begin{thm}
Any set of matrices $({\cal M}_1,\dots,{\cal M}_{n+2})$ satisfying (\ref{N6})
and such that the group $\langle{\cal M}_1,\dots,{\cal M}_{n+2}\rangle$ 
is irreducible, can be realized as monodromy group of some Fuchsian
system.
\end{thm}

\noindent{\bf Proof.} see \cite{Dek,AB}.

%%%%%%%%%%%%%%%%% REDUCTION TO GARNIER SYSTEMS

\subsection{Garnier Systems ${\cal G}_n$}

Following \cite{IKSY}, we show here how to reduce the Schlesinger 
system ${\cal S}_{n,2}$ to the Garnier system ${\cal G}_n$. 

First notice that we can always perform conformal transformations 
of the variable $\lambda$ in the Fuchsian system (\ref{N1}). For
example by $\lambda\to {\lambda-u_{n+1}\over u_{n+2}-u_{n+1}}$, 
we can fix the poles $u_{n+1}$ and $u_{n+2}$ at $0$ and $1$ 
respectively. Here, we fix them this way. 

The Schlesinger system ${\cal S}_{n,2}$ is invariant 
with respect to the gauge transformations of the form:
\begin{equation}
{\cal A}_i\sim D^{-1} {\cal A}_i D,\quad i=1,\dots,n+2,
\quad\hbox{for any $D$ diagonal matrix}.
\label{eq1}
\end{equation}
we introduce $2n$ coordinates 
$(\nu_1,\dots,\nu_{n},\rho_1,\dots,\rho_{n})$ on the quotient 
of the space of the matrices satisfying ${\cal S}_{n,2}$ with respect 
to the equivalence relation (\ref{eq1}).
The coordinates $(\nu_1,\dots,\nu_{n})$ are the roots of the 
following equation of degree $n$:
\begin{equation}
\sum_{k=1}^{n+2}{{\cal A}_{k_{12}}\over\lambda-u_k}=0, 
\label{eql}\end{equation}
and $(\rho_1,\dots,\rho_{n})$ are given by
\begin{equation}
\rho_i=\sum_{k=1}^{n+2}{{\cal A}_{k_{11}}+
{\theta_k\over2}\over\nu_i-u_k}, \qquad
i=1,\dots,n+2.
\label{eqm}\end{equation}

\begin{thm}
If $\left({\cal A}_1(u_1,\dots,u_{n}),\dots,
{\cal A}_{n+2}(u_1,\dots,u_{n})\right)$ 
satisfy ${\cal S}_{n,2}$, then the functions
$\left(\nu_1(u_1,\dots,u_{n}),
\dots,\rho_{n}(u_1,\dots,u_{n})\right)$ 
defined by (\ref{eql}), (\ref{eqm}) satisfy the Garnier system 
${\cal G}_n$.
\end{thm}

\noindent The proof of this result and the formulae expressing the
matrices ${\cal A}_1,\dots,{\cal A}_{n+2}$ in terms of the
coordinates $(\nu_1,\dots,\nu_{n},\rho_1,\dots,\rho_{n+2})$ 
can be found in \cite{IKSY}. Observe that the matrices 
${\cal A}_1,\dots,{\cal A}_{n+2}$ are classical functions of the
coordinates $(\nu_1,\dots,\nu_{n},\rho_1,\dots,\rho_{n})$.

\begin{rmk}
Observe that for $n>1$ the Garnier system ${\cal G}_n$ does not satisfy the 
Painlev\'e property. This is due to the fact that the coordinates 
$(\nu_1,\dots,\nu_{n})$ are defined as the roots of a
polynomial equation of degree $n$. It is possible to introduce a
canonical transformation 
$(\nu_1,\dots,\rho_{n})\to (q_1,\to,p_{n})$, such that the
new Hamiltonians are polynomials in $(q_1,\to,p_{n})$ and the
Painlev\'e property is satisfied (see \cite{OK}).
\end{rmk}

A different conformal transformation  
$\lambda\to {\lambda-u_{i}\over u_{j}-u_{i}}$, will correspond to a 
different choice of the position of the poles $0$ and $1$, i.e.
$u_i=0$ and $u_j=1$. The new Garnier system so obtained will be
related to the old one by a symmetry, described in the following
section.

%%%%%%%%%%%% symmetries

\subsubsection{Symmetries of the Garnier systems}
 
We give here a list (see \cite{IKSY}) of symmetries of ${\cal G}_n$, 
i.e. birational canonical transformations
$$
T:(\nu_1,\dots,\rho_{n},u_1,\dots,u_{n+2})\to
(\tilde\nu_1,\dots,\tilde\rho_{n},\tilde u_1,\dots,\tilde u_{n+2})
$$
which leave ${\cal G}_n$ invariant, modulo changes of the parameters
$$
l:(\theta_1,\dots,\theta_\infty)\to
(\tilde\theta_1,\dots,\tilde\theta_\infty).
$$
Such symmetries are easily understood as a result of a conformal
transformation of the original Fuchsian system (\ref{N1}). In the list
we choose $u_{n+1}=0$ and $u_{n+2}=1$.
$$
\begin{array}{c}
T_j\\
j=1,\dots,n\\
\end{array}
\left\{
\begin{array}{l}
\tilde\nu_i={u_j-\nu_i\over u_j-1},\\
\tilde\rho_i=-(u_j-1)\rho_i,\\
\tilde u_j={u_j\over u_j-1},\\
\tilde u_i={u_j-u_i\over u_j-1},\, i\neq j,n+1,\\
\end{array}\right.
$$
$$
T_{n+2}
\left\{
\begin{array}{l}
\tilde\nu_i=1-\nu_i,\\
\tilde\rho_i=-\rho_i,\\
\tilde u_i=1-u_i,\\
\end{array}\right.
$$
$$
T_{n+3}
\left\{
\begin{array}{l}
\tilde\nu_i={\nu_i\over\nu_i-1},\\
\tilde\rho_i=-(\nu_i-1)^2\rho_i+
{\sum_{j=1}^{n+3}\theta_j+1\over2}(\nu_i-1),\\
\tilde u_i={u_i\over u_i-1},\\
\end{array}\right.
$$
the parameters change accordingly
$$
\begin{array}{ll}
T_j:&\theta_j\leftrightarrow\theta_{n+1},\quad j=1,\dots,n,\\
T_{n+2}:&\theta_{n+2}\leftrightarrow\theta_{n+1},\\
T_{n+3}:&\theta_\infty=\theta_{n+3}\leftrightarrow\theta_{n+1}+1.\\
\end{array}
$$

%%%%%%%%%%%%% painleve' equation

\subsection{Painlev\'e VI equation}

In the case of $n=1$, the Garnier system ${\cal G}_1$ depends only on
one variable $u_1=x$, having fixed $u_2=0$ and $u_3=1$ as above. It
automatically satisfies the Painlev\'e property. In fact there is 
only one coordinate $\nu_1=y$ defined as root of a linear equation.
Let us put $\rho_1=p$, then the Garnier system ${\cal G}_1$ in 
this case is:
$$
\left\{
\begin{array}{ll}
{\partial  y\over\partial x} &=
{T( y)\over T'(x)}\left[2 p + {1\over  y-x}-
\sum_{k=1}^3{\vartheta_k\over y-u_k}\right]\\
{\partial p\over\partial x} &=-\bigg\{
T'( y) p^2 +
\big[2 y+x-\sum_ju_j-\sum_{k=1}^3
\vartheta_k(2 y+u_k-\sum_ju_j)\big]p+\\
&+{1\over4}(\sum_{k=1}^3\vartheta_k-\thi)
(\sum_{k=1}^3\vartheta_k+\thi-2)\bigg\}
{1\over T'(x)}.\\
\end{array}\right.
$$
This is the Painlev\'e sixth equation for $y(x)$ 
$$
\begin{array}{ll}
y_{xx}=&{1\over2}\left({1\over y}+{1\over y-1}+{1\over y-x}\right) y_x^2 -
\left({1\over x}+{1\over x-1}+{1\over y-x}\right)y_x+\\
&+{y(y-1)(y-x)\over x^2(x-1)^2}\left[\alpha+\beta {x\over y^2}+
\gamma{x-1\over(y-1)^2}+\delta {x(x-1)\over(y-x)^2}\right],\\
\end{array}
$$
with 
\begin{equation}
\alpha={(\thi-1)^2\over2},\quad\beta=-{\vartheta_2^2\over2},
\quad\gamma={\vartheta_3^2\over2},\quad\delta={1-\vartheta_1^2\over2}.
\label{8}\end{equation}

\vskip 0.2 cm
\begin{rmk} Observe that permutations of the poles $u_i$ 
and of the values $\vartheta_i$, $i=1,2,3,\infty$
induce transformations of $(y,x)$ of the type $x\to1-x$ and $y\to1-y$, 
$x\to{1\over x}$ and $y\to{1\over y}$, $x\to{x\over x-1}$ and
$y\to{x-y\over x-1}$ and their compositions. These
transformations are the {\it symmetries}\/ of the Painlev\`e VI
equation and correspond to the symmetries $T_3$,$T_4$ and $T_1$
of the Garnier system ${\cal G}_1$ above. Symmetries of the 
Painlev\`e VI equation have been classified by Okamoto~\cite{Ok2}. 
We list some of them in Section 3.4.  
\label{rmk0.0}\end{rmk}

\begin{rmk}
It is clear from (\ref{8}) that
changes of the signs of the parameters $\vartheta_k$, $k=1,2,3$ and
transformations on $\thi$ of type $\thi\to2-\thi$ give rise to the same
Painlev\'e VI equation.
\label{rmk0}\end{rmk}

%%%%%%%%%%%%%%%%%%%%

%%%%%%%%%%%%%%%%%%%% SECTION: Reductions of the Garnier systems

%%%%%%%%%%%%%%%%%%%%

\section{Reductions of the Garnier systems and classical solutions}

\subsection{A useful gauge transformation}

\begin{lm}
Given a Fuchsian system of the form (\ref{N1}) with residue matrices
${\cal A}_1,\dots$, ${\cal A}_{n+2}$ not all upper-triangular. 
Let $\pm{\theta_k\over 2}$ be the eigenvalues of the matrix 
${\cal A}_k$ for $k=1,\dots,n+2,\infty$. 
Given some integer $N$, if there is at least one $\theta_j\neq -2,-N,-2N$
then there exists a gauge transformation 
$G(u_1,\dots,u_{n+2},{\cal A}_1,\dots,{\cal A}_{n+2})$, rational in all 
arguments, that maps the given Fuchsian system to a new one of the same form 
(\ref{N1}) with matrices $\tilde{\cal A}_1,\dots,\tilde{\cal A}_{n+2}$
such that the new eigenvalues are $\pm{\theta_k\over 2}$ for $k\neq j$ 
and $\pm({\theta_j\over 2}+N)$. 
\label{lmb0}\end{lm}

\noindent{\bf Proof.} We give here the gauge transformation
$\Phi(\lambda)= G(\lambda)\tilde\Phi(\lambda)$ giving rise to the 
transformation $\thi\to\thi+2N$, 
$$
G_{\infty_{12}}= a,
\quad
G_{\infty_{21}}=-{\sum_{l=1}^{n+2}{\cal A}_{l_{21}}u_l^N\over\thi+N},
\quad G_{\infty_{22}}=0,
$$
and for $\thi\neq 0$
$$
G_{\infty_{11}}= \lambda^N+
{2\sum_{l=1}^{n+2}{\cal A}_{l_{22}}u_l^N\over\thi+2N}-
{(\thi+N)\sum_{l=1}^{n+2}{\cal A}_{l_{21}}u_l^{2N}\over
(\thi+2N)\sum_{l=1}^{n+2}{\cal A}_{l_{21}}u_l^N},
$$
while for $\thi=0$
$$
G_{\infty_{11}}= \lambda^N+
{\sum_{l=1}^{n+2}{\cal A}_{l_{22}}u_l^N\over N}-
{\sum_{l=1}^{n+2}{\cal A}_{l_{21}}u_l^{2N}\over
2\sum_{l=1}^{n+2}{\cal A}_{l_{21}}u_l^N}-
{\sum_{l=1}^{n+2}{\cal A}_{l_{21}}u_l^N\over 2 N^2},
$$
where $a\neq 0$ is an arbitrary parameter.
Such a gauge is always well defined in our hypotheses because
the equation $\sum_{l=1}^{n+2}{\cal A}_{l_{21}}u_l^N=0$ is
compatible with the Schlesinger equations iff 
${\cal A}_{l_{21}}=0$ for all $l=1,\dots,n+2$. 
Infact when $N=1$ and $\sum_{l=1}^{n+2}{\cal A}_{l_{21}}u_l=0$, we
obtain for all $i=1,\dots,n+2$
$$
0={\partial\over\partial u_i}\sum_{l=1}^{n+2}{\cal A}_{l_{21}}u_l=
(2+\thi) {\cal A}_{i_{21}} u_i,
$$
that for $\thi\neq-2$ implies ${\cal A}_{i_{21}}=0$ for all 
$i=1,\dots,n+2$. When $N>1$, $\sum_{l=1}^{n+2}{\cal A}_{l_{21}}u_l^N=0$
gives 
$$
0={\partial\over\partial u_i}\sum_{l=1}^{n+2}{\cal A}_{l_{21}}u_l^N=
N {\cal A}_{i_{21}} u_i^{N-1} + \sum_{s=1}^N \left[ 
\sum_{l=1}^{n+2}{\cal A}_{l}u_l^{N-s},\sum_{i=1}^{n+2} 
{\cal A}_{i} u_i^{s-1}\right]_{21},
$$
and summing on all $i$ we obtain 
$\sum_{i=1}^{n+2}{\cal A}_{i_{21}}u_l^{N-1}=0$. Iterating the same 
computation, we arrive at  $\sum_{l=1}^{n+2}{\cal A}_{l_{21}}u_l=0$,
that for $\thi\neq-2$ implies ${\cal A}_{i_{21}}=0$ for all 
$i=1,\dots,n+2$ as proved above.

The new matrices
$(\tilde{\cal A}_1,\dots,\tilde{\cal A}_{n+2})$ are given by
$\tilde{\cal A}_k=G(u_k)^{-1}{\cal A}_kG(u_k)$, $k=1,\dots,n+2$
and have parameters $\tilde\theta_k=\theta_k$, $k\neq\infty$ and
$\tilde\thi=\thi+2N$

Analogous formulae can be derived for the transformation 
$\theta_j\to\theta_j+2N$ for $j=1,\dots,n+2$. 
{\hfill$\bigtriangleup$}

\begin{lm}
Given a Fuchsian system of the form (\ref{N1}) with residue matrices
${\cal A}_1,\dots,{\cal A}_{n+2}$ not all lower-triangular. 
Let $\pm{\theta_k\over 2}$ be the eigenvalues of the matrix 
${\cal A}_k$ for $k=1,\dots,n+2,\infty$. 
Given some integer $N$, if there is at least one $\theta_j\neq 2, N, 2N,$
then there exists a gauge transformation 
$G(u_1,\dots,u_{n+2},{\cal A}_1$, $\dots$, ${\cal A}_{n+2})$, rational 
in all arguments, that maps the Fuchsian system with residue matrices
$({\cal A}_1,\dots,{\cal A}_{n+2})$ to a new one of the same form 
(\ref{N1}) with residue matrices 
$(\tilde{\cal A}_1,\dots,\tilde{\cal A}_{n+2})$
such that the new eigenvalues are $\pm{\theta_k\over 2}$ for $k\neq j$ 
and $\pm({\theta_j\over 2}-N)$. 
\label{lmb0.0}\end{lm}

\noindent{\bf Proof.} It is completely analogous to the previous 
proof.

\subsection{Smaller Monodromy Groups.}

\begin{df}
Given a Fuchsian system of the form (\ref{N1}), we say that its monodromy 
group $\langle {\cal M}_1,\dots,{\cal M}_{n+2}\rangle$ is 
{\it $l$-smaller,}\/ if $l$ monodromy matrices are equal to $\pm\ID$.
\end{df}

In this section we prove Theorem \ref{mainthmsm} that claims that
if a solution of the Garnier system ${\cal G}_{n}$ is such that 
the associated Fuchsian system has a $l$-smaller monodromy group, 
then it depends classically on the variables $u_j$ such that 
${\cal M}_j=\pm\ID$.

 %%%%%%%%%%%%% vicino a u_k

\vskip 0.1 cm
\noindent{\bf Proof of Theorem \ref{mainthmsm}.} 
First of all we want to prove that if there exists a solution 
$(\nu_1,\dots,\nu_n,\rho_1,\dots,\rho_n)$
of the Garnier system ${\cal G}_{n}$ such that the corresponding 
Fuchsian system of the form (\ref{N1}) with residue matrices 
$({\cal A}_1,\dots,{\cal A}_{n+2})$ has an $l$-smaller monodromy 
group, then there exists a solution $(\hat\nu_1,\dots$, 
$\hat\nu_{n-l}$, $\hat\rho_1$, $\dots$, $\hat\rho_{n-l})$ of the 
Garnier system ${\cal G}_{n-l}$ such that the corresponding
Fuchsian system has monodromy group generated by those matrices 
${\cal M}_1,\dots,{\cal M}_{n+2}$ that are not equal to $\pm\ID$.
We have to distinguish the two cases ${\cal M}_k=\ID$ or 
${\cal M}_k=-\ID$. We deal first with the case ${\cal M}_k=\ID$
in the additional technical hypothesis that the matrices 
$({\cal A}_1,\dots,{\cal A}_{n+2})$ are not all upper-triangular.

Suppose that one of the monodromy matrices, say ${\cal M}_k$, 
$k\neq\infty$, is equal
to the identity. This may happen if and only if the corresponding 
$\theta_k$ is an even integer. If $\theta_k= 2 K$, for some integer 
$K\neq\pm 1$, then, by Lemma \ref{lmb0}, we can map it to
$\theta_k-2(K+1)$ via a rational gauge transformation (in fact 
$\theta_k=2,2 K\neq K+1,2(K+1)$) and viceversa.  If $\theta_k=2$ we 
can simply use the constant gauge
$\left(\begin{array}{cc}0&1\\ 1&0\\ \end{array}\right)$.  
So we can assume without loss of generality that $\theta_k=-2$ for 
some $k$ . We are going to perform few gauge transformations. At each 
step we must normalise the new fundamental matrix at infinity as in
(\ref{M3}).

Let us switch the pole $u_k$ with $\infty$
by the change of variable $\lambda\to{1\over \lambda-u_k}$. To this 
aim, we have to diagonalise ${\cal A}_k$. We take
${\cal J}_k=\left(\begin{array}{cc}-1&0\\ 0&1\end{array}\right)$ 
and ${\cal G}_k$ is such that 
${\cal G}_k{\cal J}_k{\cal G}_k^{-1}={\cal A}_k$. We obtain
a new Fuchsian system with poles $0$, ${1\over u_l-u_k}$ for 
$l\neq k$, residue matrices ${\cal B}_l$ defined as follows:
$$
{\cal B}_\infty:={\cal J}_k,\qquad
{\cal B}_k:= {\cal G}_k^{-1}{\cal A}_\infty {\cal G}_k, 
\qquad\hbox{and for $l\neq k$, }\,
{\cal B}_l:= {\cal G}_k^{-1} {\cal A}_l {\cal G}_k,
$$
and monodromy matrices ${\cal N}_l$ defined as follows:
$$
{\cal N}_\infty:=\ID,\qquad
{\cal  N}_k:={\cal M}_\infty,\qquad
\hbox{and for $l\neq k$, }\,
{\cal N}_l:= {\cal M}_l.
$$
The condition ${\cal N}_\infty=\ID$ implies that the matrix 
${\cal R}_\infty$ defined in Proposition \ref{prop1} is identically 
equal to zero. This gives rise to the following equation for the 
matrix elements of ${\cal B}_l$, $l=1,\dots,n+2$
$$
\sum_{l=1}^{n+2}{\cal B}_{l_{21}} u_l^2+
\sum_{l=1}^{n+2}{\cal B}_{l_{21}} u_l\left(
\sum_{l=1}^{n+2}\left({\cal B}_{l_{22}}-
{\cal B}_{l_{11}}\right) u_l\right)\equiv 0.
$$
Keeping this relation in mind it is not difficult to prove that
the gauge transformation
$$
G(\lambda)=\left(\begin{array}{cc}
\lambda+1 & 1\\
\sum_{l=1}^{n+2}{\cal B}_{l_{21}} u_l &0\\
\end{array}\right),
$$
maps the Fuchsian system with residue matrices ${\cal B}_l$, 
$l=1,\dots n+2$, to a new Fuchsian system with residue matrices 
$\hat{\cal B}_l=G(u_l)^{-1} {\cal B}_l G(u_l)$, $l=1,\dots n+2$,
such that $\hat{\cal B}_\infty=0$ 
The Gauge $G(\lambda)$ is well defined in our hypotheses
and non-singular (because the  equation 
$\sum_{l=1}^{n+2}{\cal B}_{l_{21}} u_l=0$ is compatible with
the Schlesinger equations iff $\thi=-1$ or ${\cal B}_{l_{21}}=0$ for 
all $l=1,\dots,n+2$, both relations being false in our hypotheses).
The new Fuchsian system with matrices
$(\hat{\cal B}_1,\dots,\hat{\cal B}_{n+2})$ is not lower triangular and 
has the same monodromy matrices. 

We now want to build a solution of ${\cal G}_{n-1}$. To this aim we
switch back to a Fuchsian system with the original $u_l$, 
$l=1,\dots,n+2$. Again, we have to diagonalise 
$\hat{\cal B}_k= G(u_k)^{-1} {\cal G}_k^{-1} {\cal A}_\infty
{\cal G}_k G(u_k)$. Take $\hat{\cal G}_k$ such that 
$\hat{\cal G}_k  {\cal A}_\infty \hat{\cal G}_k^{-1} =\hat{\cal B}_k$.
The new Fuchsian system has residue matrices
$$
\hat{\cal A}_\infty={\cal A}_\infty,\qquad
\hat{\cal A}_k= 0,\qquad 
\hbox{and for $l\neq k$, }\,
\hat{\cal A}_l:= \hat{\cal G}_k^{-1} G(u_l)^{-1} {\cal G}_k^{-1} 
{\cal A}_l {\cal G}_k G(u_l) \hat{\cal G}_k ,
$$
and monodromy matrices
$$
\hat{\cal M}_\infty={\cal M}_\infty,\quad
\hat{\cal  M}_k=\ID,\quad
\hbox{and for $l\neq k$, }\,\hat{\cal M}_l= {\cal M}_l.
$$
The Fuchsian system so obtained corresponds to a solution 
$(\hat\nu_1,\dots,\hat\nu_n$, $\hat\rho_1,\dots,$ $\hat\rho_n)$ 
of the Garnier system ${\cal G}_{n-1}$ with monodromy group
$\langle{\cal M}_1,\dots,{\cal M}_{k-1},{\cal M}_{k+1}$, $\dots,$
${\cal M}_{n+2}\rangle$.

Now we want to prove that starting from the above obtained  
solution of the Garnier system ${\cal G}_{n-1}$, i.e. starting
from the Fuchsian system with residue matrices 
$(\hat{\cal A}_l,\dots,\hat{\cal A}_{k-1},\hat{\cal A}_{k+1},
\dots,\hat{\cal A}_{n+2})$, one can indeed build a one-parameter 
family of solutions
$(\tilde{\cal A}_1,\dots,\tilde{\cal A}_{n+2})$ of ${\cal G}_n$ having the
original monodromy matrices ${\cal M}_l$, $l=1,\dots,n+2$ and 
being rational functions of
$(u_1,\dots,u_{n+2},\hat{\cal A}_1,\dots,\hat{\cal A}_{n+1})$. 

Let us switch the pole $u_k$ with $\infty$ again, as above. Let
$\tilde{\cal G}_k$ be an arbitrary constant matrix, we obtain
a new Fuchsian system with poles $0$, ${1\over u_l-u_k}$ for 
$l\neq k$, residue matrices 
$$
\tilde{\cal B}_\infty=0,\qquad
\tilde{\cal B}_k= \tilde{\cal G}_k {\cal A}_\infty \tilde{\cal G}_k^{-1},
\qquad \hbox{and for $l\neq k$, }\, 
\tilde{\cal B}_l= \tilde{\cal G}_k\hat{\cal A}_l 
\tilde{\cal G}_k^{-1}.
$$
The new monodromy matrices are then
$$
\tilde{\cal N}_\infty:=\ID,\quad
\tilde{\cal N}_k:={\cal C}_k^{-1}{\cal M}_\infty{\cal C}_k,
\qquad \hbox{and for $l\neq k$, }\, 
\tilde{\cal N}_l:= {\cal C}_k^{-1}{\cal M}_l{\cal C}_k,
$$
where ${\cal C}_k$ is the connection matrix corresponding to the new 
normalisation of the fundamental matrix. If we choose 
$\tilde{\cal G}_k=\hat{\cal G}_k$ and ${\cal C}_k=\ID$ we obtain 
$\tilde{\cal B}_l=\hat{\cal B}_l$ for all $l=1,\dots,n+2,\infty$. 
Otherwise we obtain a one--parameter family of Fuchsian systems 
(it is actually a two--parameters family, but one 
parameter is factored out by diagonal conjugation) with monodromy 
group $\langle\tilde{\cal N}_1,\dots,\tilde{\cal N}_{n+2}\rangle$.

We now want to rebuild the residue matrix ${\cal B}_\infty$. This 
can be achieved with a singular gauge
$$
\tilde G(\lambda)=\left(\begin{array}{cc}
\lambda+{1\over 2}\left[\sum_{l=1}^{n+2}\tilde{\cal B}_{l_{11}} u_l
-\sum_{l=1}^{n+2}\tilde{\cal B}_{l_{22}} u_l - 
{\sum_{l=1}^{n+2}\tilde{\cal B}_{l_{12}} u_l^2\over 
\sum_{l=1}^{n+2}\tilde{\cal B}_{l_{12}} u_l} \right] &
\sum_{l=1}^{n+2}\tilde{\cal B}_{l_{12}} u_l\\
a &0\\ \end{array}\right),
$$
where $a\neq 0$ is an arbitrary parameter (factored out by diagonal 
conjugation). Such a gauge is rational in
the matrix elements of $(\tilde{\cal B}_1,\dots,\tilde{\cal B}_{n+2})$ 
and it is well defined in our hypotheses. It maps the one--parameter
family of Fuchsian systems with residue matrices $\tilde{\cal B}_l$, 
$l=1,\dots n+2$, to a new one--parameter family of Fuchsian systems 
with residue matrices 
$\check{\cal B}_l=\tilde G(u_l)^{-1}\tilde{\cal B}_l\tilde G(u_l)$, 
$l=1,\dots n+2$ and monodromy matrices 
$\check{\cal N}_l=\tilde{\cal N}_1$, $l=1,\dots n+2$.

We now want to  build to a one--parameter family of Fuchsian 
system with the original $u_l$, $l=1,\dots,n+2$. Thus we have to 
diagonalise $\check{\cal B}_k$. Let $\check{\cal G}_k$ such that
$\check{\cal B}_k\check{\cal G}_k{\cal A}_{\infty}
\check{\cal G}_k^{-1}$, we have the new residue matrices
$$
\tilde{\cal A}_{k} = 
\check{\cal G}_k^{-1}\check{\cal J}_k\check{\cal G}_k,\qquad
\tilde{\cal A}_{\infty} = {\cal A}_{\infty}
\qquad \hbox{and for $l\neq k$, }\, 
\tilde{\cal A}_{l} = 
\check{\cal G}_k^{-1}\check{\cal B}_l\check{\cal G}_k.
$$
The new monodromy matrices are then to be conjugated with the 
connection matrix corresponding to the new normalisation of the 
fundamental matrix. They coincide with the original monodromy matrices
${\cal M}_l$, $l=1,\dots,n+2$ up to diagonal conjugation. Since
$(\tilde\nu_1,\dots,\tilde\nu_n$, $\tilde\rho_1,\dots,\tilde\rho_n)$ 
are classical functions of 
$(\tilde{\cal A}_1,\dots,\tilde{\cal A}_{n+2})$ and  
$(u_1,\dots,u_{n+2})$ and 
$(\hat{\cal A}_1,\dots,\hat{\cal A}_{n+1})$ are classical
functions of 
$(\hat\nu_1,\dots,\hat\nu_n$, $\hat\rho_1,\dots,\hat\rho_n)$,
and $(u_1,\dots,u_{n+2})$, we conclude that 
$(\tilde\nu_1,\dots,\tilde\nu_n$, $\tilde\rho_1,\dots,\tilde\rho_n)$ 
are a one-parameter family of classical functions of
$(\hat\nu_1,\dots,\hat\nu_n$, $\hat\rho_1,\dots,\hat\rho_n)$,
and $(u_1,\dots,u_{n+2})$.

%%%%%%%% CASO UPPER non ho la famiglia

It is clear that when all matrices $({\cal A}_1,\dots,{\cal A}_{n+2})$ 
are upper triangular we can proceed as above, using the gauge 
transformations of the form given in Lemma \ref{lmb0.0} to reduce 
to the case of $\theta_k=2$, the same change of variables to
reduce to $\thi=2$ and a gauge transformation analogous to
$G(\lambda)$ and $\hat G(\lambda)$ to conclude. This concludes the 
proof in the case when one of the monodromy matrices is equal to 
the identity. 

Let us now suppose that one of the monodromy matrices, 
say ${\cal M}_k$, is equal to the minus identity. 
If ${\cal M}_k=-\ID$ with $k\neq\infty$, we can apply
the symmetries $T_{n+3}\cdot T_{k}$ and the change of variable
$\lambda\to {1\over\lambda-u_1}$ to map the given solution  
$({\cal A}_1,\dots,{\cal A}_{n+2})$ to a Fuchsian system with 
${\cal M}_k=\ID$. 

This concludes the proof in the case when one of the 
monodromy matrices is equal to $\pm\ID$. When more than one 
of the monodromy matrices are equal to plus or
minus the identity, we just iterate the above procedure. 
{\hfill$\bigtriangleup$}

%%%%%%%%%%%%% all'infinito

\begin{thm}
If there exists a solution $(\nu_1,\dots,\nu_n,\rho_1,\dots,\rho_n)$
of the Garnier system ${\cal G}_{n}$ such that the corresponding 
Fuchsian system of the form (\ref{N1in}) with matrices 
$({\cal A}_1,\dots,{\cal A}_{n+2})$ has monodromy matrix 
${\cal M}_\infty=\pm\ID$, then, for any $k=1,\dots,n+2$ such
that ${\cal M}_k\neq\pm\ID$, there exists a solution $(\hat\nu_1,
\dots,$ $\hat\nu_{n-1},\hat\rho_1,$ $\dots,$ $\hat\rho_{n-1})$
of the Garnier system ${\cal G}_{n-1}$ such that the corresponding
Fuchsian system  has
$\langle {\cal C}_k{\cal M}_1{\cal C}_k^{-1},\dots
{\cal C}_k{\cal M}_{k-1}{\cal C}_k^{-1},
{\cal C}_k{\cal M}_{k+1}{\cal C}_k^{-1},
\dots,{\cal C}_k{\cal M}_{n+2}{\cal C}_k^{-1}\rangle$, as monodromy group 
${\cal C}_k$ being the connection matrix of ${\cal M}_k$.
Moreover the given solution $(\nu_1,\dots,\nu_n$, $\rho_1,$ $\dots,\rho_n)$
of the Garnier system ${\cal G}_{n}$ depends classically on 
$(\hat\nu_1,\dots,\hat\nu_{n-1}$, $\hat\rho_1,
\dots,\hat\rho_{n-1}$, $u_1,\dots,u_{n+2})$. 
\label{mainthmsm1}\end{thm}

\begin{rmk}
Observe that for a solution of the Garnier system ${\cal G}_{n}$
there is always at least one monodromy matrix, say ${\cal M}_k$,
not equal to $\pm\ID$. 
\label{mdiag}\end{rmk}

\noindent{\bf Proof.} As above, we present the proof in 
the case when not all matrices ${\cal A}_l$ are upper triangular.
As seen in the proof of Theorem \ref{mainthmsm}, it is easy to
reproduce a proof for the case of all matrices ${\cal A}_l$ 
upper triangular. Suppose that ${\cal M}_\infty=\ID$ and 
${\cal M}_k\neq\pm\ID$ for all $k=1,\dots,n+2$ (if there are some 
${\cal M}_k=\pm\ID$, we first use Theorem \ref{mainthmsm}).
This condition implies that the matrix ${\cal R}_\infty$ defined 
in Proposition \ref{prop1} is identically equal to zero. This gives 
rise to the following equation for the matrix elements of 
${\cal A}_l$, $l=1,\dots,n+2$
$$
\sum_{l=1}^{n+2}{\cal A}_{l_{21}} u_l^2+
\sum_{l=1}^{n+2}{\cal A}_{l_{21}} u_l\left(
\sum_{l=1}^{n+2}\left({\cal A}_{l_{22}}-
{\cal A}_{l_{11}}\right) u_l\right)\equiv 0.
$$
Keeping this relation in mind it is not difficult to prove that
the gauge transformation
$$
G(\lambda)=\left(\begin{array}{cc}
\lambda+1 & 1\\
\sum_{l=1}^{n+2}{\cal A}_{l_{21}} u_l &0\\
\end{array}\right),
$$
maps the Fuchsian system with residue matrices ${\cal A}_l$, 
$l=1,\dots n+2$, to a new Fuchsian system with residue matrices 
$\hat{\cal A}_l=G(u_l)^{-1} {\cal A}_l G(u_l)$, $l=1,\dots n+2$
such that $\hat{\cal A}_\infty=0$. The Gauge $G(\lambda)$ is well 
defined in our hypotheses and non-singular. In fact the equation 
$\sum_{l=1}^{n+2}{\cal A}_{l_{21}} u_l=0$ is compatible with
the Schlesinger equations iff $\thi=-1$ or ${\cal A}_{l_{21}}=0$ for 
all $l=1,\dots,n+2$, both relations being false in our hypotheses.
The gauge $G(\lambda)$ preserves the monodromy matrices.

We now want to build a solution of ${\cal G}_{n-1}$. To this aim we
map one of the poles, say $u_k$ to infinity. We have 
to reduce $\hat{\cal A}_k= G(u_k)^{-1} {\cal A}_k G(u_k)$ to the
Jordan normal form ${\cal J}_k$. Let $\hat{\cal G}_k$ be such 
that $\hat{\cal A}_k=\hat{\cal G}_k{\cal J}_k\hat{\cal G}_k^{-1}$.
The new Fuchsian system has residue matrices
$$
\hat{\cal B}_\infty={\cal J}_k,\qquad
\hat{\cal B}_k= 0,\qquad 
\hbox{and for $l\neq 1$, }\,
\hat{\cal B}_l= \hat{\cal G}_k^{-1} G(u_l)^{-1} 
{\cal A}_l  G(u_l) \hat{\cal G}_k,
$$
and monodromy matrices
$$
\hat{\cal N}_\infty={\cal C}_1{\cal M}_k{\cal C}_k^{-1},\qquad
\hat{\cal  N}_k=\ID,\qquad
\hbox{and for $l\neq k$, }\,
\hat{\cal N}_l= {\cal C}_k{\cal M}_l,{\cal C}_k^{-1},
$$
where ${\cal C}_k$ is the connection matrix of ${\cal M}_k$.
Let $(\hat\nu_1,\dots,\hat\nu_{n-l},\hat\rho_1,
\dots,\hat\rho_{n-l})$ be the solution of ${\cal G}_{n-1}$
corresponding to the Fuchsian system having residue matrices
$\hat{\cal B}_l$, $l=1,\dots,k-1,k+1,\dots,n+2$.
We want to show that the given solution
$(\nu_1,\dots,\nu_n,\rho_1,\dots,\rho_n)$ of ${\cal G}_n$
depends classically on $(\hat\nu_1,\dots,\hat\nu_{n-l},\hat\rho_1,
\dots,\hat\rho_{n-l})$ and $(u_1,\dots,u_{n+2})$.

Since we can obviously recover $\hat{\cal A}_l$, $l=1,\dots,n+2$, 
from $\hat{\cal B}_l$,  $l=1,\dots,n+2$, by inverting the map above, 
we just need to show that we can reconstruct the matrices ${\cal A}_l$,
$l=1,\dots,n+2$, from $\hat{\cal A}_l$, $l=1,\dots,n+2$,  by a 
gauge transformation  that depends classically on $\hat{\cal A}_1$. 
The gauge
$$
\hat G(\lambda)=\left(\begin{array}{cc}
\lambda+{1\over 2}\left[\sum_{l=1}^{n+2}\hat{\cal A}_{l_{11}} u_l
-\sum_{l=1}^{n+2}\hat{\cal A}_{l_{22}} u_l - 
{\sum_{l=1}^{n+2}\hat{\cal A}_{l_{12}} u_l^2\over 
\sum_{l=1}^{n+2}\hat{\cal A}_{l_{12}} u_l} \right] &
\sum_{l=1}^{n+2}\hat{\cal A}_{l_{12}} u_l\\
a &0\\ \end{array}\right),
$$
where $a\neq 0$ is an arbitrary parameter (factored out by diagonal 
conjugation) is well defined and non-singular in our hypotheses. 
The new matrices $\tilde{\cal A}_l=\hat G(u_l)^{-1}\hat{\cal A}_l
\hat G(u_l)$, $l=1,\dots,n+2$, coincide with the old ones 
${\cal A}_l$, $l=1,\dots,n+2$, because of the uniqueness
Lemma (\ref{lm2.8}).

If ${\cal M}_\infty=-\ID$, by the symmetry $T_3$ and the a change of
variable $\lambda\to{1\over \lambda-u_1}$, we obtain a Fuchsian system 
with ${\cal M}_\infty=\ID$. We can then proceed as above.
{\hfill$\bigtriangleup$}

%%%%%%%%%%%%%%%%%%% monodromia riducibile

\subsection{Reducible monodromy groups}

\begin{df}
Given a Fuchsian system of the form (\ref{N1}), we say that its 
monodromy group $\langle {\cal M}_1,\dots,{\cal M}_{n+2}\rangle$ 
is {\it reducible,}\/ if the monodromy matrices admit a
common invariant subspace $X$ of dimension $1$ in the space
of solutions of the system (\ref{N1}).
\end{df}

In particular, if the monodromy group is reducible, then there
exists a basis where all monodromy matrices have the form
\begin{equation}
{\cal M}_{k}=
\left(\begin{array}{cc}
\exp(\eps_ki\pi\theta_k) & \star\\
0 & \exp(-i\eps_k\pi\theta_k)
\end{array}\right),\qquad k=1,\dots,n+2,\infty,
\label{monrid}\end{equation}
where $\eps_k=\pm 1$, the choice is fixed once for ever.

We are now going to prove Theorem \ref{mainthmred} that says that
the solutions of the Garnier systems having reducible monodromy groups
belong to $n$-parameters families of classical solutions.

\vskip 0.1 cm
\noindent{\bf Proof of Theorem \ref{mainthmred}.} First of all 
notice that by (\ref{monrid}) and (\ref{N6}), we have that 
$\thi-\sum_{k=1}^{n+2}\eps_k\theta_k\in2\interi$. 
Let $({\cal A}_1,\dots,{\cal A}_{n+2})$ be the matrices corresponding
to the given solution $(\nu_1,\dots,\nu_n,\rho_1,\dots,\rho_n)$.
We now show few Lemmata that allow us first to map the Fuchsian
system with matrices $({\cal A}_1,\dots,{\cal A}_{n+2})$ to a 
Fuchsian system with matrices upper-triangular 
$(\tilde{\cal A}_1,\dots,\tilde{\cal A}_{n+2})$ and
with the same monodromy group. We prove that the given 
$({\cal A}_1,\dots,{\cal A}_{n+2})$ are rational functions of 
$(\tilde{\cal A}_1,\dots,\tilde{\cal A}_{n+2})$ and 
$(u_1\dots,u_{n+2})$. Then we prove that the solutions
$(\tilde\nu_1,\dots,\tilde\nu_n,\tilde\rho_1,\dots,\tilde\rho_n)$
of ${\cal G}_{n}$ corresponding to the upper-triangular matrices
$(\tilde{\cal A}_1,\dots,\tilde{\cal A}_{n+2})$ are classical 
functions of $(u_1\dots,u_{n+2})$. We can then conclude 
that $(\nu_1,\dots,\nu_n,\rho_1,\dots,\rho_n)$ are classical 
functions of $(u_1\dots,u_{n+2})$. We follow this procedure first in the 
case when no monodromy matrix equal to $\pm\ID$, and finally we allow
$l$--smaller monodromy groups too.

\begin{lm}
If there exists a solution $(\nu_1,\dots,\nu_n,\rho_1,\dots,\rho_n)$
of the Garnier system ${\cal G}_{n}$ such that the corresponding 
Fuchsian system has a reducible monodromy group with no monodromy
matrix equal to $\pm\ID$, and 
$\thi-\sum_{k=1}^{n+2}\eps_k\theta_k=2 K$ for some integer
$K\neq0$, then there exists a solution
$(\hat\nu_1,\dots,\hat\nu_n,$ $\hat\rho_1,\dots,\hat\rho_n)$
of ${\cal G}_{n}$ with parameters 
$\hat{\theta}_\infty-\sum_{k=1}^{n+2}\eps_k\hat\theta_k=0$
and with the same monodromy data. The matrices
$({\cal A}_1,\dots,{\cal A}_{n+2})$ are rational functions of
$(\hat{\cal A}_1,\dots,\hat{\cal A}_{n+2})$ and $(u_1\dots,u_{n+2})$.
\label{lmb2}\end{lm}

\noindent{\bf Proof.} 
Suppose $\thi-\sum_{k=1}^{n+2}\eps_k\theta_k\in 2K\neq0$, then 
the residue matrices $({\cal A}_1,\dots,{\cal A}_{n+2})$ are not all 
upper-triangular (lower-triangular) otherwise $K=0$. Suppose now 
that there are at least two parameters, say $\theta_i$, $\theta_j$ 
that are non-zero and not equal to $2$. Given $K$, we can find two 
integers $N_i$, $N_j$ 
such that $\eps_i N_i+\eps_j N_j=-K$ and 
$\theta_{i,j}\neq -N_{i,j},-2 N_{i,j}, -2 (N_{i,j}+1)$. As a 
consequence the gauge transformations 
$G_i$ and $G_j$, $i\neq j$, defined in Lemma \ref{lmb0} are well defined. 
In fact, chosen any two matrices ${\cal A}_i$, ${\cal A}_j$, 
$i\neq j$, the number of integers for which the gauge transformations 
$G_i(\lambda)$ and $G_j(\lambda)$ are not defined is $4$, but there 
are infinite integers $N_i$, $N_j$ such that $\eps_i N_i+\eps_j N_j=-K$.
Then applying $G_i(\lambda)$ and $G_j(\lambda)$, we find a new
Fuchsian system
with matrices $(\hat{\cal A}_1,\dots,\hat{\cal A}_{n+2})$ 
such that $\hat\theta_k=\theta_k$ for $k\neq i,j$,  
$\hat\theta_j=\theta_j+2N_j$ and $\hat\theta_j=\theta_j+2N_i$. 
Then $\hat{\theta}_\infty-\sum_{k=1}^{n+2}\eps_k\hat\theta_k\in2\interi=0$,
as we wanted. The obtained matrices
$(\hat{\cal A}_1,\dots,\hat{\cal A}_{n+2})$ are not all lower triangular 
because otherwise the monodromy group would be lower triangular as 
well. Moreover $\hat\theta_{i,j}\neq 2, N_{i,j}, 2 N_{i,j}$ because 
of the choice of $N_{i,j}$.
Then we can apply Lemma \ref{lmb0.0} and build two gauges 
$\hat G_i(\lambda)$ and $\hat G_j(\lambda)$, that depend rationally on 
$(\hat{\cal A}_1,\dots,\hat{\cal A}_{n+2})$. The new matrices
$\tilde{\cal A}_l=G_i(u_l)^{-1}G_j(u_l)^{-1}\hat{\cal A}_l
G_j(u_l)G_i(u_l)$, $l=1,\dots,n+2$ must coincide with the original 
ones $({\cal A}_1,\dots,{\cal A}_{n+2})$ because they have the same 
monodromy matrices the same eigenvalues and the same poles (see Lemma 
\ref{lm2.8}). So the matrices
$({\cal A}_1,\dots,{\cal A}_{n+2})$ are rational functions of
$(\hat{\cal A}_1,\dots,\hat{\cal A}_{n+2})$ and $u_1\dots,u_{n+2}$.

If there are not two parameters $\theta_i$ $\theta_j$ that are 
non-zero, there must be at
least one non zero, say $\theta_i$. In fact if all $\theta_k$ were 
zero then $K=0$. So suppose $\theta_i=2 K$. We apply 
$G_j(\lambda)$ defined in Lemma \ref{lmb0} and obtain 
$\hat\theta_j=0\neq 2,K,2K$. Again we can apply Lemma \ref{lmb0.0} 
and build another gauge $\hat G_i(\lambda)$ that depend rationally on 
$(\hat{\cal A}_1,\dots,\hat{\cal A}_{n+2})$ and allows to reconstruct
the original residue matrices $({\cal A}_1,\dots,{\cal A}_{n+2})$. 

Suppose that all $\theta_k=-2$. We
apply the gauge $\left(\begin{array}{cc}0&1\\ 1&0\\
\end{array}\right)$ to invert the sign of $\theta_j$ and then 
proceed as above. {\hfill$\bigtriangleup$}

\begin{lm}
Let the the Fuchsian system (\ref{N1}) have an reducible monodromy
group and let the sum $\thi-\sum_{k=1}^{n+2}\eps_k\theta_k$ be zero. 
Then there exists a gauge transformation $\Phi=P\tilde\Phi$ 
independent on $\lambda,u_1,\dots,u_{n+2}$ such that
$$
{{\rm d}\over{\rm d}\lambda}\tilde\Phi=
\sum_{k=1}^{n+2}{\tilde{\cal A}_k\over \lambda-u_k}\tilde\Phi,
$$
has the upper-triangular form
$$
\tilde{\cal A}_{k_{21}}=0,\quad\hbox{for all }
k=1,\dots,n+2.
$$
\label{lmb3}\end{lm}
\noindent{\bf Proof.} see \cite{AB}.

\begin{lm}
The solutions 
$(\tilde\nu_1,\dots,\tilde\nu_n,\tilde\rho_1,\dots,\tilde\rho_n)$
of the Garnier system ${\cal G}_{n}$ with corresponding Fuchsian
system of the upper-triangular form
$$
\tilde{\cal A}_{k_{21}}=0,\quad\hbox{for all }
k=1,\dots,n+2.
$$
are classical solutions and can be expressed via Lauricella Hypergeometric 
equations.
\label{lmlauricella}\end{lm}

\noindent{\bf Proof.} Suppose that $\tilde{\cal A}_{k_{21}}=0$ for all 
$k=1,\dots,n+2$. Then we can choose 
$\tilde{\cal A}_{k_{21}}=-{\theta_k\over2}$, and thus the
corresponding coordinates $(\tilde\rho_1,\dots,\tilde\rho_n)$
are all identically equal to zero. As shown in \cite{OK}, the 
coordinates $(\tilde\nu_1,\dots,\tilde\nu_n)$ are then classical 
functions and can be expressed via Lauricella Hypergeometric 
equations. {\hfill$\bigtriangleup$}

Putting together Lemmata \ref{lmb2}, \ref{lmb3}, we have shown that if there 
exists a solution $(\nu_1,\dots,\nu_n,\rho_1,\dots,\rho_n)$ of the 
Garnier system ${\cal G}_{n}$ such that the corresponding Fuchsian 
system has a reducible monodromy group with no monodromy matrix equal 
to $\pm\ID$, and $\thi-\sum_{k=1}^{n+2}\eps_k\theta_k=2 K$ for some 
integer $K\neq0$, then there exists a solution
$(\tilde\nu_1,\dots,\tilde\nu_n,$ $\tilde\rho_1,\dots,\tilde\rho_n)$
of ${\cal G}_{n}$ with parameters 
$\hat\theta_1,\dots,\hat\theta_{n+2},\hat{\theta}_\infty$ given in the proof
of Lemma \ref{lmb2},
with upper triangular residue matrices and with the same monodromy
data. The solution $(\nu_1,\dots,\nu_n,\rho_1,\dots,\rho_n)$ of 
the Garnier system ${\cal G}_{n}$ is then classical by Lemma 
\ref{lmlauricella}. We want to show that this solution belongs
to a $n$-parameter family of classical solutions of ${\cal G}_{n}$.

Given the parameters $\hat\theta_1,\dots,\hat\theta_{n+2},
\hat{\theta}_\infty$ such that 
$\hat{\theta}_\infty-\sum_{k=1}^{n+2}\eps_k\hat\theta_k=0$, we can
build a $(n+1)$-parameter family of upper-triangular Fuchsian systems with
residue matrices
$$
{\cal A}_l=\left(\begin{array}{cc}{\eps_l\theta_l\over2}&0\\
0&-{\eps_l\theta_l\over2}\end{array}\right).
$$
Applying the gauges $\hat G_i(\lambda)$ and $\hat G_j(\lambda)$ 
defined in Lemma \ref{lmb0.0} we obtain a $(n+1)$-parameter family
of Fuchsian systems with the initial exponents 
$\theta_1,\dots,\theta_{n+2},\thi$ and reducible monodromy group. 
This leads to the existence of a
$n$-parameter family (one parameter is factored out by diagonal 
conjugation) of classical solutions of the Garnier system ${\cal G}_{n}$.

This concludes the proof of Theorem \ref{mainthmred} in the case of
no monodromy matrix equal to $\pm\ID$. In the case of $l$-smaller 
reducible monodromy groups we have the following 

\begin{lm}
If there exists a solution $(\nu_1,\dots,\nu_n,\rho_1,\dots,\rho_n)$
of the Garnier system ${\cal G}_{n}$ such that the corresponding 
Fuchsian system has an $l$-smaller reducible monodromy group  
then it is a classical solution.
\label{cor1}\end{lm}

\noindent{\bf Proof.} 
Suppose that the Fuchsian system corresponding to the solution
$(\nu_1,\dots,\nu_n,\rho_1,\dots,\rho_n)$ is non upper or lower
triangular (otherwise we immediately obtain that the solution
$(\nu_1,\dots,\nu_n,\rho_1,\dots,\rho_n)$ is a classical solution 
by Lemma \ref{lmlauricella}). Then we can apply Theorem
\ref{mainthmsm}. Then there exists a solution
$(\tilde\nu_1,\dots,\tilde\nu_{n-l}$, 
$\tilde\rho_1,\dots,\tilde\rho_{n-l})$ of ${\cal G}_{n-l}$ 
with reducible monodromy group and
$(\nu_1,\dots,\nu_n,\rho_1,\dots,\rho_n)$ are classical 
functions of 
$(\tilde\nu_1,\dots,\tilde\nu_{n-l}$, 
$\tilde\rho_1,\dots,\tilde\rho_{n-l})$ and $u_1,\dots,u_{n+2})$. 
{\hfill$\bigtriangleup$}

\subsection{Classical Solutions of the Painlev\'e VI equation}

In the case of the Painlev\'e VI equation some symmetries
other than $T_{1,3,4}$ are known, namely the transformations 
$w_i$ which preserve $x$:
$$
\begin{array}{cl}
w_0(y,p,b_1, b_2, b_3, b_4)&=
(y,p-{b_3+b_4+1\over y-x},b_1,b_2,-1-b_4,-1-b_3)\\
w_1(y,p,b_1, b_2, b_3, b_4)&=(y,p-{b_1-b_2\over y-1},b_2, b_1, b_3, b_4),\\
w_2(y,p,b_1, b_2, b_3, b_4)&=(\tilde y,\tilde p,b_1, b_3, b_2, b_4),\\
w_3(y,p,b_1, b_2, b_3, b_4)&=(y,p,b_1, b_2, b_4, b_3),\\
w_4(y,p,b_1, b_2, b_3, b_4)&=(y,p-{b_1+b_2\over y},-b_2,-b_1,b_3,b_4),\\
\end{array}
$$
where, in the case when $u_1=0$, $u_2=x$ ad $u_3=1$,
$$
b_1={\theta_1+\theta_3\over 2},\quad
b_2={\theta_1-\theta_3\over 2},\quad
b_3={\theta_2+\theta_\infty\over 2}-1,\quad
b_4={\theta_2-\theta_\infty\over 2}.
$$ 
The formulae for $\tilde y(y,p,x)$ and $\tilde p(y,p,x)$ are given in 
\cite{Ok2}. The above transformations $w_0,\dots,w_4$ generate a group 
$\tilde W$ that is isomorphic to $W_a(D_4)$ the affine extension of
the Weyl group of $D_4$. 

Let $R$ be the collection of the following $24$ vectors\footnote{Root
system of type $D_4$.}
$(\pm1,\pm1,0,0)$, $(\pm1,0,\pm1,0)$, $(\pm1,0,0,\pm1)$,
$(0,\pm1,\pm1,0)$, $(0,\pm1,0,\pm1)$, $(0,0,\pm1,\pm1)$. 
Consider the hyper-planes in the space of the parameters 
${\bf b}=(b_1,b_2, b_3, b_4)\in\complessi^4$ 
$$
H_{r,k}=\{{\bf b}\in\complessi^4|\, ({\bf b},{\bf r})=k\},
$$
where ${\bf r}\in R$ and $k\in\interi$. Define the following sets
$$
\begin{array}{c}
M= \bigcup_{k\in\interi,\, r\in R}H_{r,k},\\
P= \bigcup_{k_1,k_2\in\interi,\, r_1,r_2\in R}
H_{r_1,k_1}\cap H_{r_2,k_2},\\
L= \bigcup_{k_1,k_2,k_3\in\interi,\, r_1,r_2,r_3\in R}
H_{r_1,k_1}\cap H_{r_2,k_2}\cap H_{r_3,k_3},\\
D=\bigcup_{k_1,k_2,k_3,k_4\in\interi,\, r_1,r_2,r_3,r_4\in R}
H_{r_1,k_1}\cap H_{r_2,k_2}\cap H_{r_3,k_3}\cap H_{r_4,k_4},
\end{array}
$$
where the unions over $2$ or more elements in $R$ are to be considered
only among linearly independent elements of $R$. Obviously 
$D\subset L\subset P\subset M$. In terms of the parameters 
$(\theta_1,\theta_2,\theta_3,\thi)$ the sets $H_{r,k}$ can be
rewritten as 
$$
\tilde H_{i,k}=\left\{(\theta_1,\theta_2,\theta_3,\thi)
\in\complessi^4|\, \theta_i=k\right\}, \qquad i=1,2,3,\infty,
$$
$$
\tilde H_{0,k}=\left\{(\theta_1,\theta_2,\theta_3,\thi)
\in\complessi^4|\,\sum_{j=1}^3\eps_j\theta_j+\thi=k\right\}, 
\qquad \eps_j=\pm,1
$$
so that the sets $D\subset L\subset P\subset M$ have the form
$$
\begin{array}{c}
M= \bigcup_{k\in\interi,\, i=0,1,2,3,\infty}\tilde H_{i,k},\\
P= \bigcup_{k_1,k_2\in\interi,\,\, i_1,i_2=0,1,2,3,\infty}
\tilde H_{r_1,k_1}\cap\tilde H_{r_2,k_2},\\
L= \bigcup_{k_1,k_2,k_3\in\interi,\,i_1,i_2,i_3 =0,1,2,3,\infty}
\tilde H_{r_1,k_1}\cap\tilde H_{r_2,k_2}\cap\tilde H_{r_3,k_3},\\
D=\bigcup_{k_1,k_2,k_3,k_4\in\interi,\,
\,i_1,i_2,i_3,i_4 =0,1,2,3,\infty}
\tilde H_{r_1,k_1}\cap\tilde H_{r_2,k_2}\cap\tilde H_{r_3,k_3}
\cap\tilde H_{r_4,k_4},
\end{array}
$$
where $i_l\neq i_m$ for $l\neq m$.

Classical non-algebraic solutions of the Painlev\'e VI equation are 
classified in the following (see \cite{Wat})

\begin{thm}
\begin{enumerate}
\item For ${\bf b}\in M\setminus P$ there exists a one-parameter
family of classical solutions of the Painlev\'e VI equation.
\item For ${\bf b}\in P\setminus L$ there exist two one-parameter
families of classical solutions of the Painlev\'e VI equation.
\item For ${\bf b}\in L\setminus D$ there exist three one-parameter
families of classical solutions of the Painlev\'e VI equation.
\item For ${\bf b}\in D$ there exist four one-parameter families 
of classical solutions of the Painlev\'e VI equation,
\item Any non-algebraic solution of the Painlev\'e VI 
equation with ${\bf b}\not\in M$ is non-classical.
\end{enumerate}
\label{thw1}\end{thm}

Using Theorem \ref{thw1} it is possible to prove Theorem
\ref{thmio}. The proof is based on two Lemmata that classify
all solutions of the Painlev\'e VI equation having reducible or 
$1$-smaller monodromy groups. 

\begin{lm} All solutions of the Painlev\'e VI corresponding to reducible
monodromy groups are equivalent via birational transformations
generated by $w_0,w_1,w_3,w_4$ and via symmetries $T_1$,$T_3$ and $T_4$
to the following one-parameter family of solutions, realized for 
$\thi=-(\theta_1+\theta_2+\theta_3)$
\begin{equation}
y={(1+\theta_1+\theta_2-x-\theta_2 x) u(x) -x(x-1) u_x(x)
\over (1+\theta_1+\theta_2+\theta_3) u(x)}\label{ric}
\end{equation}
where $u(x)=u_1(x)+\nu u_2(x)$, $u_1(x),u_2(x)$ being two linear independent
solutions of the following hypergeometric equation 
\begin{equation}
\begin{array}{cl} x(x-1) u_{xx}(x)=&\left[2+\theta_1+\theta_2-
(4+\theta_1+2\theta_2+\theta_3)x\right] u_x(x)-\\
&-(2+\theta_1+\theta_2+\theta_3)(\theta_2+1) u(x).\\ 
\end{array} \label{hyp}
\end{equation}
\label{lmred}\end{lm}

\noindent{\bf Proof.}
For reducible monodromy groups there exists a basis in which all
monodromy matrices are upper triangular. We can always perform
a change of basis in order that ${\cal M}_\infty$ has the form (\ref{N9}) 
and all monodromy matrices have the form:
$$
{\cal M}_k=\left(\begin{array}{cc}
\exp(\pi i\theta_k)&\star\\  0&\exp(-\pi i\theta_k)\\ 
\end{array}\right).
$$ 
It then follows, by the relation (\ref{N6}), that 
$\thi+\varepsilon_k\sum_k\theta_k=2 N$, $\varepsilon_k=\pm1$, $N\in\interi$.
By means of the birational transformations $w_0,w_1,w_4$,
which trivially preserve $y(x)$, we can always assume that $\varepsilon_k=+1$.
Perform the following gauge transformation on the Fuchsian system
$$
\Phi\to\tilde\Phi=\Pi_{k=1}^3(\lambda-u_k)^{-\theta_k\over2}\Phi,
\qquad
{\cal A}_k\to\tilde{\cal A}_k={\cal A}_k-{\theta_k\over2}\ID.
$$ 
The new residue at infinity is, for $\thi\neq0$,
$$
\tilde{\cal A}_\infty=\left(\begin{array}{cc}
{\thi+\sum\theta_k\over 2}&0\\
0&{-\thi+\sum\theta_k\over 2}\\ \end{array}\right),
$$
while for $\thi=0$, is
$$
\tilde{\cal A}_\infty=\left(\begin{array}{cc}
{\sum\theta_k\over 2}&1\\
0&{\sum\theta_k\over 2}\\ \end{array}\right).
$$
Then the new monodromy matrices are
$$
\tilde{\cal M}_k=\exp(-\pi i\theta_k){\cal M}_k
=\left(\begin{array}{cc}1&\star\\ 0&\exp(-2\pi i\theta_k)\\
\end{array}\right),\qquad k=1,2,3,
$$
and 
$$
\tilde{\cal M}_\infty=\left(\begin{array}{cc}1&0\\ 0& 
\exp[\pi i(-\thi+\sum_k\theta_k)]\\
\end{array}\right)\exp({\cal R}_\infty).
$$
Due to the form of the monodromy matrices,
the new Fuchsian system admits a non-zero
single valued solution $\tilde Y$. This solution is analytic at $u_1$,
$u_2$ and $u_3$ because all residue matrices $\tilde{\cal A}_k$, $k=1,2,3$
have a zero eigenvalue. At infinity we have polynomial behaviour.
Applying the birational transformations generated by $w_0,w_1,w_3,w_4$
and the symmetries $T_1$,$T_3$ and $T_4$, we can assume $N=0$ without
loss of generality. For $N=0$, $\tilde Y$ is necessarily constant. 
Thus, near each $u_j$, one has
$$
{\tilde{\cal A}_j\over\lambda-u_j} \tilde Y+{\cal O}(1)=0,
$$
that implies that $\tilde Y$ has the form 
$\left(\begin{array}{c}a\\ 0\\ \end{array}\right)$, for
some $a\neq0$ and all residue matrices ${\cal A}_k$ are upper
triangular. Correspondingly $p=\sum_{k=1}^3{\theta_k\over(q-u_k)}$ 
and the solution $y(x)$ of the Painlev\'e VI is a solution of the 
following Riccati equation
$$
y_x(x)={1+\theta_1+\theta_2+\theta_3\over x(x-1)} y^2-
{1+\theta_1+\theta_2+\theta_1 x+\theta_3 x\over x(x-1)} y+
{\theta_1\over x-1}\label{ric1}
$$
which can be solved by means of the hypergeometric equation
$$
y={(1+\theta_1+\theta_2-x-\theta_2 x)  -x(x-1) {u_x(x)\over u(x)}
\over 1+\theta_1+\theta_2+\theta_3}
$$
where $u(x)$ is a solution of the hypergeometric equation with 
$$
a=1+\theta_2,\quad b=2+\theta_1+\theta_2+\theta_3, 
\quad c=2+\theta_1+\theta_2.\label{ric2}
$$
This concludes the proof of the theorem.
{\hfill$\bigtriangleup$}

\begin{lm} 
All solutions of the Painlev\'e VI corresponding to $1$--smaller
monodromy groups belong to the following list
\begin{enumerate}
\item {\it Forbidden solutions}\/ $y(x)\equiv 0,x,1$, realized for 
$\theta_1=0$, $\theta_2=0$, $\theta_3=0$ respectively. The corresponding
monodromy data are such that ${\cal M}_1=\ID$, ${\cal M}_2=\ID$, 
${\cal M}_3=\ID$, respectively. 
\item {\it Forbidden solution}\/ $y(x)\equiv\infty$ realized for
$\thi=1$. The corresponding monodromy data are such that
$M_\infty=-\ID$.
\item {\it Generalised Chazy solutions.}\/ They are realized in
each of the following four cases:  
$\thi$ odd integer, $\thi\neq 1$ and ${\cal M}_\infty=-\ID$ or
$\theta_k$ even integer, $\theta_k\neq0$ and ${\cal M}_k=\ID$, for each
$k=1,2,3$. For $\thi=-1$, they have the following form
\begin{equation}
\begin{array}{cl}
y(x)=&-x \left[(W+2 x-\theta_1 W-\theta_3 W)^2-4 x+W(4
\theta_3-\theta_2^2 W)\right]\cdot\\ 
&\left\{(x-1)\left[(W+2 x-\theta_1 W-\theta_3 W)^2-4 x+W(4
\theta_3-\theta_2^2 W)\right]-4 \theta_1 W^2\right\}\cdot\\ 
&\big\{[(\theta_1+\theta_3-1)^2-\theta_2^2][\theta_3(x-1)+\theta_1 x]W^3
-2[(x-1)(3x-2)\theta_3^2+\\ 
&+2 (x-1)(1 - 2 x)\theta_3+x(3x-1)\theta_1^2
+2 x(1-2x)\theta_1- \\ 
&-x(x-1)\theta_2^2+6x(x-1)\theta_1\theta_3] W^2+\\
&+4x(x-1)[3\theta_3(x-1)+3\theta_1 x+1-2x] W-8 x^2(x-1)^2\big\}^{-1}
\\ \end{array}
\label{DH1}\end{equation}
where $W={u(x)\over u'(x)}$, $u(x)$ being is any solution of the
following hypergeometric equation:
\begin{equation}
x(1-x)u''+ [1-\theta_3-(2-\theta_1-\theta_3)x]u'-
{(\theta_1+\theta_3-1)^2-\theta_2^2\over4} u=0.
\label{DH2}\end{equation}
\item {\it Riccati-type solutions.}\/ They are all solutions of the 
Painlev\'e VI equation realized for $\theta_k=2 n+1$, $n\in\interi$, 
for one or more $k=1,2,3$, or for $\theta_\infty=2 n$, $n\in\interi$, 
that can be obtained by
the birational transformations generated by $w_0,w_1,w_3,w_4$ and 
the symmetries $T_1$,$T_3$ and $T_4$ from the following one-parameter 
family of solutions of the Painlev\'e VI realized for $\thi=2$
\begin{equation}
y(x)=f(x){f(x)[-\theta_1+\theta_2+\theta_3-x(\theta_1+\theta_2+\theta_3)]
+2\theta_1 x\over (-\theta_1+\theta_2-\theta_3) f(x)^2 +
2(\theta_3-\theta_2 x)f(x)+ (\theta_1+\theta_2-\theta_3)x}
\label{rt1}\end{equation}
where $f(x)$ is a solution of the following Riccati equation
\begin{equation}
f_x={f(f-1)(f-x)\over x(x-1)}
\left[{\theta_3-\theta_2-\theta_1\over 2 f}+
{\theta_1+\theta_2+\theta_3+2\over 2(f-x)}+
{\theta_1-\theta_2-\theta_3\over 2(f-1)}\right]. 
\label{rt}\end{equation}
They are realized for
monodromy matrices are $M_\infty=\ID$ and $M_k=-\ID$ respectively.
\end{enumerate}
\label{lmchazy}\end{lm}

\begin{rmk}
In the limit $\theta_k\to0$ for all
$k=1,2,3$, the above generalised Chazy solutions become the Chazy
solutions found in [Ma].
\label{rmkcha}\end{rmk}

\noindent{\bf Proof.} Let us consider the case 
${\cal M}_\infty=\pm\ID$ (the case ${\cal M}_k=\pm\ID$ for
some $k=1,2,3$ can be obtained from this by means of the symmetries
$T_1$,$T_3$ and $T_4$). First consider the case ${\cal M}_\infty=\ID$. 
Then the corresponding $\thi$ is an even integer. Applying the birational 
transformations generated by $w_0,w_1,w_3,w_4$ and the symmetries 
$T_1$,$T_3$ and $T_4$, we can choose $\theta_\infty=2$ without
loss of generality. The assumption ${\cal M}_\infty=\ID$ implies that 
${\cal R}_{\infty_{12}}=0$, i.e., for $\thi=2$,
$$
T(q) p^2-P(q) p 
\sum_{k=1}^3{\theta_k\over q-u_k} +{\theta_1+\theta_2+\theta_3\over 4} 
\sum_{k=1}^3\theta_k (q+u_k-u_i-u_j)=0.
$$
By substitution we obtain an algebraic relation 
${\cal R}_\infty(x,y(x),p(x))=0$ between $y(x)$ and its conjugate
momentum $p(x)$, whose common solutions with the Painlev\'e VI
equation are all given by (\ref{rt1}) and (\ref{rt}).
Let ${\cal M}_\infty=-\ID$. Then $\thi\in 2\interi+1$ and 
${\cal R}_\infty=0$. The equation ${\cal R}_\infty=0$ for $\thi=1$
leads to ${\cal A}_{k_{12}} u_k=0$,
which gives the forbidden solution $y=\infty$. All other cases with
even integer $\thi\neq 1$ are equivalent via birational
transformations to the case $\thi=-1$. 
For $\thi=-1$, the equation ${\cal R}_\infty=0$ corresponds to
${\cal A}_{k_{21}} u_k=0$, i.e. to the following algebraic relation
between $y(x)$ and $p(x)$:
$$
b_0(y,x)p^4+b_1(y,x) p^3+b_2(y,x)p^2+b_3(y,x)p+b_4(y,x)=0,
$$
where 
$$
\begin{array}{cl}
b_0 =& 16 P(y)^4,\\ 
b_1 = &-32 P(y)\left[\theta_1(y-1)(y-x)+\theta_2 y(y-1)+
\theta_3 y (y-x)\right]\\ 
b_2 =&8 \big[-P(y)(3y-1-x)\thi^2+(y-1)(y-x)(3 y^2-3 y+2x-3y x)
\theta_1^2+\\ 
&+y(y-1)(3 y^2-3 y+x-x^2)\theta_2^2+y(y-x)
(3 y^2-1+x-3 y x)\theta_3^2+\\ 
&+6 P(y)\left[\theta_1 \theta_2 (y-1)+
\theta_1 \theta_3(y-x)+y \theta_2 \theta_3\right]\big]\\ 
\end{array}
$$
$$
\begin{array}{cl}
b_3=&-8\bigg\{2 P(y)\thi^3 - (3y-1-x)
\big[\theta_1 (y-1)(y-x)+\theta_2 y(y-1)+\\
&+\theta_3y(y-x)\big]\thi^2+
\left[\theta_1 (y-1-x)+\theta_2(y+x-1)+ \theta_3(y+1-x)\right]\cdot\\ 
&\cdot(\theta_1+\theta_2+\theta_3)
\left[\theta_1(y-1)(y-x)+ \theta_2y(y-1) +\theta_3y(y-x)\right]
\bigg\}\\ \end{array}
$$
$$
\begin{array}{cl}
b_4=&(\theta_1+\theta_2+\theta_3-\thi)\big\{
(3y^2-2y-1+2x-2yx-x^2)\thi^3+\\ 
&+\big[(6y-1-5y^2-6x+6yx-x^2)\theta_1+
(6y-1-5y^2+\\ 
&+2x-2yx-x^2)\theta_2+(2x-2y-1-5y^2+6yx-x^2)\theta_3\big]\thi^2+\\ 
&+(\thi+\theta_1+\theta_2+\theta_3)
\left[\theta_1 (y-1-x)+\theta_2(y+x-1)+ \theta_3(y+1-x)
\right]^2\big\}.
\\ \end{array}
$$
Such a relation gives an algebraic differential equation of the first
order that is satisfied only by the generalised Chazy solutions.
{\hfill$\bigtriangleup$}

\vskip 0.1 cm
\noindent{\bf Proof of Theorem \ref{thmio}.} 

Observe that the birational transformations generated by 
$w_0,w_1,w_2,w_3,w_4$ and the symmetries $T_1$,$T_3$ and $T_4$ 
preserve each of the spaces $M$, $P$, $L$ and $D$. Such symmetries 
preserve the class of the classical solutions. In the Lemmata 
\ref{lmchazy} and \ref{lmred} we have found that every time
one $\theta_k\in\interi$ or 
$\thi+\sum_{k=1}^3\eps_k\theta_k\in 2\interi$ there exists a 
classical solution of the Painlev\'e VI. Thanks to Theorem
\ref{thw1}, the solutions of Lemmata \ref{lmchazy} and \ref{lmred} give all 
one parameter families of classical solutions classified in Theorem
\ref{thw1}.
{\hfill$\bigtriangleup$}

\bibliography{garnier.bib}
\bibliographystyle{plain}

\end{document}